%% file: main.tex
\documentclass[10pt]{article}

\usepackage{definitions}
\title{PDNQP: A GPU-based Factorization-free Method for Large-scale Nonconvex Quadratic Programming}
\author{Zixi Chen\thanks{Peking University, School of Mathematical Sciences (chenzixi22@pku.edu.cn).} \and Haihao Lu\thanks{MIT, Sloan School of Management (haihao@mit.edu).}}
\date{}

\begin{document}
\maketitle
\begin{abstract}
Large-scale nonconvex quadratic programming remains challenging because the sparse matrix factorizations used by many existing solvers can incur substantial computational and memory costs and are difficult to parallelize efficiently on GPUs. We present \PDNQP{}, a factorization-free first-order solver for finding stationary points of large-scale nonconvex quadratic programs. The method combines a proximal augmented-Lagrangian framework with restarted accelerated primal--dual hybrid gradient to solve the resulting strongly convex QP subproblems. Its key ingredients are a residual-form reformulation that avoids penalty-dependent normal matrices and preserves a fixed sparse constraint operator across outer iterations, and an adaptive strategy that coordinates inner-solve accuracy with outer progress. The resulting computations rely on sparse matrix--vector products, projections, and vector operations suited to GPU execution. On $107$ nonconvex CUTEst instances, \PDNQP{} achieves the highest success rate among the tested solvers at both accuracy levels. With absolute and relative tolerances set to $10^{-4}$, it achieves speedups of $2.6\times$--$6.0\times$ in shifted geometric mean runtime; at $10^{-6}$, its aggregate runtime remains comparable to that of the fastest tested solvers. On $30$ large-scale synthetic nonconvex QPs with approximately one million variables each, \PDNQP{} returns points satisfying the common external termination criteria on $27$ instances, while none of the other tested solvers succeeds under the same experimental protocol.
\end{abstract}

\section{Introduction}
\label{sec:introduction}

Quadratic programming (QP) is a fundamental class of optimization problems with applications across operations research, machine learning, control, and scientific computing. In this paper, we consider linearly constrained QPs with a possibly indefinite quadratic objective. Such problems arise in graph-based soft clustering \cite{yu2005softclustering}, network flow models with concave costs \cite{zangwill1968concavecost}, and support vector machines with indefinite kernels \cite{liu2016indefinitekernel}. Indefinite QPs also arise as subproblems in sequential quadratic programming when the Hessian of the Lagrangian is not positive semidefinite \cite{boggs1995sqp}, and in economic model predictive control with indefinite quadratic stage costs \cite{zanon2014indefinitempc}. Moreover, classical continuous formulations of the maximum-clique problem \cite{motzkin1965maxima} and Nash equilibria in bimatrix games \cite{mangasarian1964games} take the form of nonconvex QPs.

Solving nonconvex QPs at scale remains challenging. While global optimization methods are available for small problems, large-scale solvers typically seek stationary points. Practical algorithms for this task commonly rely on Newton-type or other second-order steps that require solving sparse linear systems. For example, Ipopt implements a primal--dual interior-point method that approximately solves a sequence of barrier subproblems using Newton steps computed via sparse matrix factorizations \cite{wachter2006ipopt}. Another particularly relevant example is QPALM, a proximal augmented-Lagrangian method that solves a sequence of strongly convex QP subproblems \cite{hermans2022qpalm}. QPALM computes semismooth Newton directions for these subproblems using sparse matrix factorizations and factorization updates. Both approaches rely on sparse matrix factorizations, which can incur substantial memory and computational costs for very large problems and can be difficult to parallelize efficiently on GPUs.

Recent progress in large-scale convex optimization suggests a different approach. First-order methods primarily use sparse matrix--vector products and simple vector operations, avoiding the storage and computation of matrix factorizations and aligning naturally with modern parallel hardware. PDLP demonstrated that first-order methods, when equipped with suitable preconditioning, adaptive step sizes, and restarting, can be effective for large-scale linear programming \cite{applegate2021pdlp}. cuPDLP subsequently showed that this factorization-free computational structure can be efficiently implemented on GPUs \cite{lu2025cupdlp,lu2025cupdlpx}. PDQP extended this approach to convex QP using restarted accelerated primal--dual hybrid gradient (rAPDHG) \cite{lu2025practicalqp}. More recent methods, including PDHCG and HPR-QP, further demonstrate the scalability of first-order and matrix-free methods for large-scale convex QPs on GPUs \cite{huang2025pdhcg,chen2025hprqp}. Related developments extend this direction to broader convex quadratic and conic problems \cite{lin2025pdcs,li2026pdhcgii,li2026pdhcgcqp}.

The goal of this paper is to extend this scalable factorization-free approach to nonconvex QP. We adopt a proximal augmented-Lagrangian framework in which a sufficiently strong proximal regularization makes each inner QP strongly convex. These inner problems can then, in principle, be solved by a factorization-free first-order method. One key challenge is to formulate the inner QPs in a way that preserves the sparse structure of the original problem and enables its efficient reuse across outer iterations. A second challenge is to adapt inner-solve accuracy to outer progress, avoiding unnecessarily accurate inner solves while maintaining effective outer updates.

A compact formulation of the proximal augmented-Lagrangian subproblem eliminates the constraint residuals and retains only the primal and slack variables. Although this reduces the dimension of the inner problem, it introduces penalty-weighted normal operators such as $\sigma A^\top A$ into the primal Hessian. Forming these matrices explicitly can produce substantial fill-in, even when the original constraint matrix is sparse. Moreover, the resulting normal operator changes whenever the outer penalty parameters are updated, limiting the reuse of a fixed sparse operator across inner solves.
Our central idea is to retain the constraint residuals as explicit variables. This increases the dimension of the inner QPs but avoids introducing penalty-dependent normal matrices into the primal Hessian. More importantly, the resulting strongly convex QPs share the same sparse constraint operator across outer iterations. They can therefore be solved by rAPDHG using only sparse matrix--vector products, projections, and elementary vector operations.

To address the second challenge, \PDNQP{} uses an adaptive inner stopping rule that coordinates subproblem accuracy with outer progress. Since high-accuracy subproblem solves can be costly for a first-order method, the inner tolerances are initially loose and tightened gradually. The rate of tightening depends on the magnitude of the proximal residual relative to the outer projected-stationarity residual: the tolerances decrease slowly when this ratio is small and more aggressively as it grows. This strategy avoids excessive inner work in early outer iterations while providing increasingly accurate subproblem solutions as the outer method approaches stationarity.

Building on the above, we develop \PDNQP{}, a factorization-free first-order solver for large-scale nonconvex QP. The outer method follows a proximal augmented-Lagrangian framework, while the strongly convex inner QPs are solved by rAPDHG. We incorporate diagonal scaling, adaptive outer and inner updates, warm starts, restarting, and Anderson acceleration to improve practical performance. The computationally intensive inner iterations are executed on the GPU and reuse the same sparse operator structure throughout the outer iterations. Provided that the inner-solve errors decrease sufficiently fast, the base method retains the convergence guarantees of the proximal augmented-Lagrangian framework \cite{hermans2022qpalm}.

We evaluate \PDNQP{} on two collections of nonconvex QPs. The first consists of $107$ nonconvex instances from CUTEst \cite{gould2015cutest}. The second consists of $30$ application-inspired instances with approximately one million variables each, derived from graph soft-labeling, concave-cost multistage flow, and indefinite graph-kernel SVM models. We compare \PDNQP{} with QPALM \cite{hermans2022qpalm}, Ipopt \cite{wachter2006ipopt}, MadNLPGPU \cite{shin2024gpu}, and Knitro \cite{byrd2006knitro}.
At a relative tolerance of $10^{-4}$, \PDNQP{} successfully solves all $107$ CUTEst instances; at a relative tolerance of $10^{-6}$, it solves $101$ instances. On the million-variable test set, \PDNQP{} returns externally verified solutions for $27$ of the $30$ instances, whereas none of the other four solvers returns a verified solution. These results indicate that \PDNQP{} combines robustness on standard nonconvex benchmarks with the ability to solve QPs at a scale beyond the reach of the other tested methods.

The main contributions of this paper are as follows.
\begin{itemize}
    \item \textbf{Residual-form proximal augmented-Lagrangian framework.} We introduce a residual-form formulation in which the successive strongly convex inner QPs share a fixed sparse constraint operator. This formulation avoids penalty-dependent normal matrices and provides the structural foundation for factorization-free computation.

    \item \textbf{Factorization-free GPU solver.}
We develop \PDNQP{}, which solves the inner QPs using rAPDHG and incorporates diagonal scaling, adaptive parameter updates, warm starts, restarting, and Anderson acceleration. Its dominant computations consist of sparse matrix--vector products, projections, and vector operations that are well suited to GPU execution.

\item \textbf{Large-scale numerical evaluation.}
We evaluate \PDNQP{} on standard nonconvex CUTEst instances and on application-inspired QPs with approximately one million variables. The results demonstrate strong performance on standard benchmarks and scalability to problem sizes that the other tested solvers do not successfully handle.
\end{itemize}

\subsection{Related literature}
\label{subsec:related-literature}

\paragraph{Local optimization of nonconvex QP.}
A substantial literature studies local solutions or stationary points of nonconvex QPs.
Early affine-scaling methods were developed by Ye for nonconvex quadratic programming
\cite{ye1992affinescaling}.
Ye also studied the complexity of computing approximate KKT points of general QPs
\cite{ye1998kkt}.
Tseng analyzed Dikin-type affine-scaling methods for nonconvex quadratic minimization
\cite{tseng2004dikin}, and Absil and Tits developed Newton--KKT interior-point methods for possibly indefinite QPs
\cite{absil2007newtonkkt}.
These methods establish an important algorithmic foundation for nonconvex QP, but their main computational steps are based on Newton systems, affine-scaling linear algebra, or related second-order computations.

\paragraph{Global optimization of nonconvex QP.}
A separate line of work seeks globally optimal solutions.
Representative approaches include reformulation--convexification methods
\cite{sherali1995reformulation},
branch-and-bound methods based on semidefinite relaxations
\cite{burer2008branch},
and reformulations based on completely positive programming
\cite{chen2012completelypositive}.
These methods provide global guarantees but generally require substantially more computation than methods that seek stationary points.
Our objective in this paper is the latter: we target stationary points of large-scale nonconvex QPs rather than globally optimal solutions.

\paragraph{Proximal augmented-Lagrangian methods.}
The connection between augmented-Lagrangian and proximal-point methods goes back to Rockafellar
\cite{rockafellar1976augmented}.
Proximal regularization is particularly useful for nonconvex QPs because it can make the optimization subproblem strongly convex.
QPALM develops this idea into a practical solver for nonconvex QPs
\cite{hermans2022qpalm}.
It applies a proximal augmented-Lagrangian framework and solves the resulting strongly convex subproblems using semismooth Newton directions and sparse factorization updates.
The method establishes convergence to stationary points under suitable conditions and provides an effective approach for nonconvex QPs of moderate size.
\PDNQP{} adopts the same general proximal augmented-Lagrangian principle, but uses a different representation and a factorization-free first-order method for the inner subproblems.

\paragraph{Traditional convex quadratic programming methods.}
Classical algorithms for convex QP include active-set, interior-point, augmented-Lagrangian, and operator-splitting methods.
Active-set and interior-point methods typically obtain their search directions by solving linear systems, while operator-splitting methods use simpler iterations but may still require a factorization or an inner linear solve.
For example, OSQP applies ADMM to convex QPs and repeatedly solves a quasi-definite linear system whose factorization can be reused
\cite{stellato2020osqp}.
Methods for convex and more general QPs have also been developed using active-set strategies
\cite{gill2015generalqp}.
These approaches are effective for a wide range of QPs, but the cost of sparse linear-system solves can become important as the problem size grows.

\paragraph{Large-scale first-order and GPU methods.}
Recent work on LP and convex QP has shown that first-order methods can provide an alternative computational structure for very large problems.
PDLP develops restarted PDHG into a practical large-scale LP solver using diagonal preconditioning, adaptive step sizes, and restarting
\cite{applegate2021pdlp}.
cuPDLP and its subsequent implementations exploit the matrix--vector structure of these iterations to obtain efficient GPU solvers
\cite{lu2025cupdlp,lu2025cupdlpx}.
PDQP extends this direction to convex QP and develops restarted accelerated PDHG for quadratic objectives
\cite{lu2025practicalqp}.

Several recent methods further develop scalable first-order algorithms for convex QP.
PDHCG combines a primal--dual framework with inexact conjugate-gradient solution of the primal subproblems and provides an efficient GPU implementation
\cite{huang2025pdhcg}.
HPR-QP develops a dual Halpern Peaceman--Rachford method for large-scale convex composite QPs and is also implemented on GPUs
\cite{chen2025hprqp}.
PDCS considers large-scale conic optimization, including convex QPs, using restarted primal--dual iterations and GPU-oriented sparse matrix--vector computation
\cite{lin2025pdcs}.
More recently, PDHCG-II improves PDHCG using Halpern-type acceleration, adaptive primal--dual weighting, and adaptive inner stopping
\cite{li2026pdhcgii}, while PDHCG-CQP extends this line to large-scale conic convex QPs and multi-GPU computation
\cite{li2026pdhcgcqp}.
Together, these works demonstrate that matrix-free and first-order methods can scale effectively on modern parallel hardware.
However, they address convex quadratic objectives, whereas the matrix $Q$ in the present work may be indefinite.

\paragraph{General-purpose nonlinear programming solvers.}
Nonconvex QPs can also be solved by general-purpose nonlinear programming solvers.
Ipopt uses a primal--dual interior-point method
\cite{wachter2006ipopt}, and Knitro provides several nonlinear optimization algorithms, including interior-point methods
\cite{byrd2006knitro}.
MadNLP and MadNLPGPU use an interior-point framework with GPU-accelerated linear algebra
\cite{shin2024gpu}.
These solvers are designed for much broader nonlinear optimization models and rely on linear-system solves within their second-order frameworks.

\subsection{Notations and paper organization}
The symbol $I$ denotes an identity matrix of the required dimension.  For a vector $a$, the positive part $(a)_+$ is defined componentwise by $((a)_+)_i=\max\{a_i,0\}$.  For a box $[\ell,u]$, $\proj_{[\ell,u]}$ denotes the Euclidean projection.  Let $\mathcal D\subseteq\R^n$ be a nonempty closed convex set.  For $x\in\mathcal D$, its normal cone is
\[
 N_{\mathcal D}(x)
 =\{v\in\R^n:\langle v,y-x\rangle\le0\ \text{for every }y\in\mathcal D\}.
\]
We set $N_{\mathcal D}(x)=\varnothing$ when $x\notin\mathcal D$.  

\paragraph{Paper organization.}
\Cref{sec:base} introduces the nonconvex QP formulation and presents the base \PDNQP{} method.
\Cref{sec:practical} describes the practical enhancements and GPU implementation.
\Cref{sec:experiments} presents the numerical experiments.

\section{Problem formulation and the base \PDNQP{} method}
\label{sec:base}

We consider the nonconvex QP of the following form:
\begin{equation}
\label{eq:qp}
\begin{aligned}
    \min_{x\in\R^n}\quad &
    f(x) \coloneqq c^\top x+\tfrac{1}{2} x^\top Qx,\\
    \text{s.t.}\quad & A_{\E}x=b_{\E},\\
    & A_{\I}x\ge b_{\I},\\
    & \ell\le x\le u \ ,
\end{aligned}
\end{equation}
where $n$ is the number of variables, and $m_{\E}$ and $m_{\I}$ are the numbers of equality and inequality constraints, respectively. The problem data consist of a symmetric, possibly indefinite matrix $Q\in\R^{n\times n}$; vectors $c\in\R^n$, $b_{\E}\in\R^{m_{\E}}$, and $b_{\I}\in\R^{m_{\I}}$; matrices $A_{\E}\in\R^{m_{\E}\times n}$ and $A_{\I}\in\R^{m_{\I}\times n}$; and bound vectors $\ell,u\in(\R\cup{-\infty,+\infty})^n$ satisfying $\ell\le u$.  
The equality and inequality row sets are denoted by $\E$ and $\I$, respectively.

Let $\lambda_{\E}\in\R^{m_{\E}}$ and $\lambda_{\I}\in\R^{m_{\I}}$ denote the equality and inequality multipliers, respectively. We aim to find a stationary point (i.e., a KKT point) $(x,\lambda_{\E},\lambda_{\I})$ satisfying the stationarity condition
\begin{equation}
\label{eq:kkt-stationarity}
0\in Qx+c+A_{\E}^{\top}\lambda_{\E}
+A_{\I}^{\top}\lambda_{\I}+N_{[\ell,u]}(x),
\end{equation}
together with the primal feasibility constraints in \eqref{eq:qp}, dual feasibility $\lambda_{\I}\le0$, and complementarity $\lambda_{\I}^{\top}(A_{\I}x-b_{\I})=0$. Here, $N_{[\ell,u]}(x)$ denotes the normal cone to the box $[\ell,u]$ at $x$.

\PDNQP{} seeks such stationary points by combining a proximal augmented-Lagrangian outer framework with the factorization-free restarted accelerated primal--dual hybrid gradient (rAPDHG) method of \cite{lu2025practicalqp} for the inner solves. A formulation that eliminates the constraint residuals introduces penalty-weighted normal operators involving $A_{\E}^{\top}A_{\E}$ and $A_{\I}^{\top}A_{\I}$. Explicitly forming these operators can produce substantial fill-in, and their numerical values change whenever the outer penalty parameters are updated.

The central algorithmic idea in \PDNQP{} is to retain the constraint residuals as explicit variables.  The resulting inner problem is a convex QP whose sparse constraint operator remains fixed across outer iterations.  Penalty updates affect only diagonal blocks of its Hessian, whereas center and multiplier updates affect only its linear term.  This structure permits the inner problems to be solved using sparse matrix--vector products and projections, with efficient reuse of sparse data structures on the GPU.

This section develops the base method of \PDNQP{} in four steps.  We first introduce the proximal ALM subproblem, then derive its residual form, describe the inner rAPDHG iteration, and finally assemble the complete outer iteration.  The adaptive rules and implementation enhancements used in the practical solver are deferred to \cref{sec:practical}.

\subsection{Proximal augmented-Lagrangian subproblem}
\label{subsec:base-outer}

Introduce a slack variable $s_{\I}\in\R^{m_{\I}}$ and rewrite the inequality constraints as $s_{\I}\ge b_{\I}$ and $A_{\I}x-s_{\I}=0$.  At outer iteration $k$, let $\lambda_{\E}^k\in\R^{m_{\E}}$ and $\lambda_{\I}^k\in\R^{m_{\I}}$ be the current outer dual variables.  Their concatenation is $\lambda^k=(\lambda_{\E}^k,\lambda_{\I}^k)\in\R^{m_{\E}+m_{\I}}$.  Let $\sigma^k>0$ be the scalar penalty.  The augmented Lagrangian is
\begin{equation}
\label{eq:base-alm}
\begin{aligned}
 \mathcal L_{\sigma^k}(x,s_{\I};\lambda^k)
 ={}&\tfrac{1}{2}x^\top Qx+c^\top x
 +(\lambda_{\E}^k)^\top(A_{\E}x-b_{\E})
 +\tfrac{\sigma^k}{2}\normtwo{A_{\E}x-b_{\E}}^2 +(\lambda_{\I}^k)^\top(A_{\I}x-s_{\I})
 +\tfrac{\sigma^k}{2}\normtwo{A_{\I}x-s_{\I}}^2.
\end{aligned}
\end{equation}
A classical ALM iteration minimizes \eqref{eq:base-alm} over
$\ell\le x\le u$ and $s_{\I}\ge b_{\I}$, and then updates the
multipliers using the resulting constraint residuals.

Because $Q$ may be indefinite, this ALM subproblem need not be convex.
To obtain a strongly convex subproblem, proximal ALM adds a quadratic
regularization term centered at $z^k\in\R^n$. Let $\gamma>0$ satisfy
\begin{equation}
\label{eq:base-strong-convexity}
 Q+\gamma^{-1}I\succ0.
\end{equation}
Using the center $z^k$, the dual variables $\lambda^k$, the penalty $\sigma^k$, and the step $\gamma$, the method approximately solves
\begin{equation}
\label{eq:base-proximal-alm}
\begin{aligned}
 \min_{\substack{\ell\le x\le u\\s_{\I}\ge b_{\I}}}\quad
 &\mathcal L_{\sigma^k}(x,s_{\I};\lambda^k)
 +\frac{1}{2\gamma}\normtwo{x-z^k}^2.
\end{aligned}
\end{equation}
Condition \eqref{eq:base-strong-convexity} and $\sigma^k>0$ make \eqref{eq:base-proximal-alm} strongly convex.  The proximal term serves two purposes: it convexifies the subproblem and limits the displacement from the current center.

Let $x^{k+1}$ be the primal point returned by the inner solver.
Minimizing \eqref{eq:base-proximal-alm} with respect to $s_{\I}$ gives
\[
s_{\I}^{k+1}
=
\max\left\{
b_{\I},
A_{\I}x^{k+1}+\frac{\lambda_{\I}^k}{\sigma^k}
\right\},
\]
where the maximum is taken componentwise. The multiplier updates are
therefore
\begin{subequations}
\label{eq:projected-multipliers}
\begin{align}
\lambda_{\E}^{k+1}
&=\lambda_{\E}^k
+\sigma^k(A_{\E}x^{k+1}-b_{\E}),\\
\lambda_{\I}^{k+1}
&=\min\left\{
0,\,
\lambda_{\I}^k
+\sigma^k(A_{\I}x^{k+1}-b_{\I})
\right\}.
\end{align}
\end{subequations}
Thus, the equality multipliers remain unrestricted, whereas the
componentwise projection in the inequality update preserves
$\lambda_{\I}^{k+1}\le0$.

If \eqref{eq:base-proximal-alm} is solved exactly, its first-order
optimality condition with respect to $x$ can be written as
\[
0\in
Qx^{k+1}+c
+A_{\E}^{\top}\lambda_{\E}^{k+1}
+A_{\I}^{\top}\lambda_{\I}^{k+1}
+\gamma^{-1}(x^{k+1}-z^k)
+N_{[\ell,u]}(x^{k+1}).
\]
This condition differs from the stationarity condition of the original
QP only by the proximal residual
$\gamma^{-1}(x^{k+1}-z^k)$. Consequently, if
$x^{k+1}-z^k\to0$, stationarity of the proximal subproblems yields
stationarity of the original QP.

\subsection{Residual-form convex subproblem}
\label{subsec:residual-form}

Let $r_{\E}\in\R^{m_{\E}}$ be the equality residual variable.  Let $r_{\I}\in\R^{m_{\I}}$ be the inequality-split residual variable.  Eliminating these variables would add the following matrix to the $x$-block of the Hessian:
\[
 \sigma^k\bigl(A_{\E}^{\top}A_{\E}+A_{\I}^{\top}A_{\I}\bigr).
\]
Even when $A_{\E}$ and $A_{\I}$
are sparse, this normal matrix can be substantially denser. Moreover,
it changes whenever the penalty parameter $\sigma^k$ is updated.
To avoid forming this matrix, we retain $r_{\E}$ and $r_{\I}$ as
explicit variables and impose their defining equations as constraints.
This gives the equivalent subproblem
\begin{equation}
\label{eq:prox-alm-subproblem}
\begin{aligned}
 \min_{x,s_{\I},r_{\E},r_{\I}}\quad
 &\tfrac{1}{2}x^\top Qx+c^\top x
 +\frac{1}{2\gamma}\normtwo{x-z^k}^2\\
 &+(\lambda_{\E}^k)^\top r_{\E}
 +\tfrac{\sigma^k}{2}\normtwo{r_{\E}}^2
 +(\lambda_{\I}^k)^\top r_{\I}
 +\tfrac{\sigma^k}{2}\normtwo{r_{\I}}^2,\\
 \text{s.t.}\quad
 &A_{\E}x-r_{\E}=b_{\E},\\
 &A_{\I}x-s_{\I}-r_{\I}=0,\\
 &\ell\le x\le u,\qquad s_{\I}\ge b_{\I}.
\end{aligned}
\end{equation}
Eliminating $r_{\E}$ and $r_{\I}$ recovers \eqref{eq:base-proximal-alm}.  Thus the residual form changes only the representation of the subproblem, not the proximal ALM step.

To express \eqref{eq:prox-alm-subproblem} as a standard convex QP, set $d_w=n+m_{\E}+2m_{\I}$.  Define the inner primal vector
\[
 w=(x,s_{\I},r_{\E},r_{\I})\in\R^{d_w}.
\]
Let $H_k\in\R^{d_w\times d_w}$ be the inner Hessian.  Let $q_k\in\R^{d_w}$ be the inner linear term.  Let $B\in\R^{(m_{\E}+m_{\I})\times d_w}$ be the inner constraint matrix, and let $h\in\R^{m_{\E}+m_{\I}}$ be its right-hand side.  Let $\underline w,\overline w\in(\R\cup\{-\infty,+\infty\})^{d_w}$ be the lower and upper bound vectors.  The subproblem can be written compactly as
\begin{equation}
\label{eq:inner-qp}
 \min_{\underline w\le w\le\overline w}
 \quad
 \tfrac{1}{2}w^\top H_kw+q_k^\top w
 \qquad
 \text{s.t.}\quad Bw=h.
\end{equation}
The matrix $B$ is fixed, and $h$ is its right-hand side:
\[
 B=
 \begin{bmatrix}
 A_{\E}&0&-I&0\\
 A_{\I}&-I&0&-I
 \end{bmatrix},
 \qquad
 h=
 \begin{bmatrix}
 b_{\E}\\0
 \end{bmatrix}.
\]
The Hessian $H_k$ and linear term $q_k$ are
\[
 H_k=
 \diag\!\left(
 Q+\gamma^{-1}I,\,0,\,\sigma^k I,\,\sigma^k I
 \right),
 \qquad
 q_k=
 \begin{bmatrix}
 c-\gamma^{-1}z^k\\
 0\\
 \lambda_{\E}^k\\
 \lambda_{\I}^k
 \end{bmatrix}.
\]
The componentwise lower and upper bounds are
\begin{equation}
\label{eq:inner-bounds}
 \underline w=(\ell,b_{\I},-\infty,-\infty),
 \qquad
 \overline w=(u,+\infty,+\infty,+\infty).
\end{equation}
The constant $(2\gamma)^{-1}\normtwo{z^k}^2$ has been omitted because
it does not affect the minimizer.  Under \eqref{eq:base-strong-convexity}, $H_k\succeq0$ and is positive definite on $\ker B$.  Hence \eqref{eq:inner-qp} has a unique minimizer.

The key benefit of this formulation is that the constraint operator
$B$ and right-hand side $h$ are independent of the outer iteration
$k$. A penalty update changes only the diagonal residual blocks of
$H_k$, while updates of the proximal center and outer multipliers
change only $q_k$. Therefore, every outer iteration produces a convex
QP with the same sparse constraint operator. This fixed-operator
structure enables the inner problems to be solved without
factorization, using only sparse matrix--vector products and
componentwise projections, while reusing the same GPU sparse-matrix
data structures throughout the solve.

\subsection{Inner restarted accelerated PDHG}
\label{subsec:base-pdhg}

Let $\mathcal W=[\underline w,\overline w]\subseteq\R^{d_w}$.  Let $y\in\R^{m_{\E}+m_{\I}}$ be the inner dual variable for $Bw=h$.  Let $\mathcal L_k:\mathcal W\times\R^{m_{\E}+m_{\I}}\to\R$ be the saddle function.  The saddle-point formulation of \eqref{eq:inner-qp} becomes
\begin{equation}
\label{eq:base-inner-saddle}
 \min_{w\in\mathcal W}\max_{y\in\R^{m_{\E}+m_{\I}}}
 \mathcal L_k(w,y),
 \qquad
 \mathcal L_k(w,y)
 =
 \tfrac{1}{2}w^\top H_kw+q_k^\top w-y^\top(Bw-h).
\end{equation}
The saddle problem \eqref{eq:base-inner-saddle} is convex--concave, and its only nonsmooth primal term is the box indicator with a componentwise projection.  It is therefore natural to use a primal--dual first-order method.  We apply the rAPDHG method of \cite{lu2025practicalqp}: acceleration exploits the quadratic objective, while restarting periodically replaces the long running average by a new reference point.

The inner dual variable $y$ is different from the outer dual variables $\lambda_{\E}^k$ and $\lambda_{\I}^k$.  For a fixed outer subproblem, let $j$ index restart epochs and let $t$ index iterations in one epoch.  The vectors $w^{j,t},\bar w^{j,t}\in\R^{d_w}$ are the current and averaged inner primal iterates.  The vectors $y^{j,t},\bar y^{j,t}\in\R^{m_{\E}+m_{\I}}$ are the current and averaged inner dual iterates.  The vector $w_{\mathrm{md}}^{j,t}\in\R^{d_w}$ is the middle point used for the quadratic gradient.  The parameters $\beta_t>0$ and $\theta_t\ge0$ are scalar averaging and extrapolation parameters.  The scalars $\tau_t>0$ and $\eta_t>0$ are the primal and dual step sizes.  The symbol $\gamma$ remains reserved for the proximal step.

\begin{algorithm}[H]
\caption{Restarted accelerated PDHG for \eqref{eq:inner-qp}}
\label{alg:base-pdhg}
\begin{algorithmic}[1]
\Require Initial vectors $w^{0,0}\in\R^{d_w}$ and $y^{0,0}\in\R^{m_{\E}+m_{\I}}$; scalar parameters $\{(\beta_t,\theta_t,\tau_t,\eta_t)\}_{t=0}^{K-1}$; restart length $K$
\State $j\gets0$
\Repeat
    \State
    $(\bar w^{j,0},\bar y^{j,0})
    \gets(w^{j,0},y^{j,0})$,
    $w^{j,-1}\gets w^{j,0}$
    \For{$t=0,\ldots,K-1$}
        \State
        $
        w_{\mathrm{md}}^{j,t}
        \gets
        (1-\beta_t^{-1})\bar w^{j,t}
        +\beta_t^{-1}w^{j,t}
        $
        \State
        $
        y^{j,t+1}
        \gets
        y^{j,t}
        +\eta_t
        \left[
        h-B\left(
        w^{j,t}
        +\theta_t(w^{j,t}-w^{j,t-1})
        \right)
        \right]
        $
        \State
        $
        w^{j,t+1}
        \gets
        \proj_{\mathcal W}\!\left(
        w^{j,t}
        -\tau_t\left(
        H_kw_{\mathrm{md}}^{j,t}
        +q_k-B^\top y^{j,t+1}
        \right)
        \right)
        $
        \State
        $
        \bar w^{j,t+1}
        \gets
        (1-\beta_t^{-1})\bar w^{j,t}
        +\beta_t^{-1}w^{j,t+1}
        $
        \State
        $
        \bar y^{j,t+1}
        \gets
        (1-\beta_t^{-1})\bar y^{j,t}
        +\beta_t^{-1}y^{j,t+1}
        $
    \EndFor
    \State
    $(w^{j+1,0},y^{j+1,0})
    \gets(\bar w^{j,K},\bar y^{j,K})$
    \State $j\gets j+1$
\Until{$(w^{j,0},y^{j,0})$ satisfies the prescribed inner stopping test}
\State \Return $(w^{j,0},y^{j,0})$
\end{algorithmic}
\end{algorithm}

\Cref{alg:base-pdhg} has two levels.  Within a restart epoch, the middle point supplies the accelerated quadratic gradient, while the running average stabilizes the primal--dual iterates.  After $K$ iterations, the averaged point becomes the initial point of the next epoch.  Each inner iteration consists only of products with $H_k$, $B$, and $B^\top$, together with a componentwise projection onto $\mathcal W$.

The practical inner solver in \cref{subsec:practical-inner} preserves this two-level structure.  It replaces the scalar metric and fixed restart frequency by diagonal adaptive steps, an acceptance test, adaptive restarting, and state reuse across outer subproblems.

\subsection{Complete base iteration}
\label{subsec:complete-base}

We now combine the proximal ALM outer step with the residual-form rAPDHG inner solve.  At outer iteration $k$, the current center, multipliers, and penalty define one convex QP of the form \eqref{eq:inner-qp}.  After this QP has been solved to the prescribed accuracy, its $x$-block is used to update the outer multipliers.  The outer rule then determines whether to update the proximal center and how to choose the next penalty and inner tolerance. Algorithm \ref{alg:base-outer} presents this base \PDNQP{} algorithm. It leaves the center-acceptance, penalty-update, accuracy-update rules, and other practical choices in \cref{sec:practical}.

\begin{algorithm}[H]
\caption{Base \PDNQP{} method}
\label{alg:base-outer}
\begin{algorithmic}[1]
\Require Initial center $z^0\in\R^n$; outer dual variables $\lambda_{\E}^0\in\R^{m_{\E}}$ and $\lambda_{\I}^0\in\R^{m_{\I}}$ with $\lambda_{\I}^0\le0$; scalar penalty $\sigma^0>0$; scalar proximal step $\gamma$ satisfying \eqref{eq:base-strong-convexity}; prescribed inner and outer accuracies
\For{$k=0,1,\ldots$}
    \State Form the convex subproblem \eqref{eq:inner-qp}
    \State Approximately solve \eqref{eq:inner-qp} by
    \cref{alg:base-pdhg} to the prescribed inner accuracy and obtain $w^{k+1}$ with $x$-block $x^{k+1}$
    \State Update $(\lambda_{\E}^{k+1},\lambda_{\I}^{k+1})$ by \eqref{eq:projected-multipliers}
    \If{the requested accuracy for \eqref{eq:qp} is satisfied}
        \State \Return $(x^{k+1},\lambda_{\E}^{k+1},\lambda_{\I}^{k+1})$
    \EndIf
    \State Choose $z^{k+1}$, $\sigma^{k+1}>0$, and the next inner accuracy according to the outer update rule
\EndFor
\end{algorithmic}
\end{algorithm}

\paragraph{Convergence guarantees.}
The residual-form reformulation changes only the representation of
each proximal ALM subproblem. Indeed, eliminating $r_{\E}$ and
$r_{\I}$ from \eqref{eq:prox-alm-subproblem} recovers exactly the
proximal ALM subproblem \eqref{eq:base-proximal-alm}. Consequently,
provided that the inner solutions satisfy the same inexactness
conditions, the outer iteration of the base \PDNQP{} method falls
within the proximal ALM framework analyzed for QPALM
\cite{hermans2022qpalm}.
In particular, suppose that the original QP is feasible and its
objective is bounded below on the feasible set. Assume further that
\[
Q+\gamma^{-1}I\succ0,
\qquad
0<\underline{\sigma}\le \sigma^k\le\overline{\sigma}<+\infty,
\]
and that the inner-solve errors satisfy the inexactness conditions in
\cite{hermans2022qpalm}. The convergence results established there
then apply directly to the base \PDNQP{} iteration. In particular, if
the inner errors decrease geometrically, the outer iterates converge
to a stationary point of \eqref{eq:qp}, and the accepted primal
iterates converge $R$-linearly. We therefore omit a separate proof.

\section{Practical enhancements}
\label{sec:practical}

The base PDNQP method presented in the previous section specifies the proximal augmented-Lagrangian framework and the rAPDHG inner iteration. To improve the practical performance, this section presents a series of enhancements that balance the numerical scales of the problem, coordinate progress between the outer and inner iterations, and reuse the common structure of successive inner QPs.   These enhancements are important for the numerical performance of PDNQP. 

\subsection{Initialization and diagonal scaling}
\label{subsec:practical-setup}

The initialization phase is performed once at the beginning of the algorithm.  It balances the objective and constraints, selects a proximal step that convexifies the inner QPs, and initializes the primal variable, outer dual variables, row penalties, and stopping tolerances used by the practical method.

\subsubsection{Outer-loop initialization}
The initial primal variable is $x^0\in\R^n$, and the initial outer dual variables are $\lambda_{\E}^0\in\R^{m_{\E}}$ and $\lambda_{\I}^0\in\R^{m_{\I}}$.  User-provided values are used when available; otherwise, these vectors are set to zero.  The initial proximal center is set to $z^0=x^0$.

The outer loop of PDNQP uses row-wise penalties instead of one scalar penalty as described in Algorithm \ref{alg:base-outer}. The equality and inequality penalty vectors $\sigma_{\E}^0$ and $\sigma_{\I}^0$ are initialized using the same method in QPALM \cite{hermans2022qpalm}.  This selection aims to balance the infeasibility of each coordinate.

The scalars $\delta_{a,k}$ and $\delta_{r,k}$ are the inner stopping tolerances, while $\varepsilon_{a,k}$ and $\varepsilon_{r,k}$ are the center-update tolerances.  They are initialized with $\delta_{a,0}=\delta_{r,0}=\varepsilon_{a,0}=\varepsilon_{r,0}=1$ and are updated during the outer iterations.

\subsubsection{Selection of the proximal regularization ($\gamma$)}

The proximal regularization term $\gamma$ is chosen so that the objective of the ALM subproblem~\eqref{eq:prox-alm-subproblem} is convex. Specifically, it is chosen as

\begin{equation}
\label{eq:R-choice}
 \gamma=\frac{1}{-\widehat\lambda_Q+
       \delta_\gamma\max\{1,\norminf{Q}\}},
\end{equation}
where $\widehat\lambda_Q$ is an estimation of the smallest eigenvalue of the matrix $Q$ by Lanczos iterations and $\delta_{\gamma}=10^{-2}$ is a hyperparameter of the implementation.

\subsubsection{Diagonal scaling of the objective and constraints}
First-order methods are sensitive to the condition number of the optimization problem.  \PDNQP{} applies $8$ iterations of Ruiz scaling to the combined constraint matrix \cite{ruiz2001scaling} to better balance the variables and constraints.  Set
\[
 A=
 \begin{bmatrix}
 A_{\E}\\A_{\I}
 \end{bmatrix},
\]
and let $D_c$ and $D_x$ be the resulting positive diagonal row and column scaling matrices from Ruiz scaling.  The scaled constraint matrix is
\begin{equation}
\label{eq:scaled-problem}
 \bar A=D_cAD_x.
\end{equation}
The same row and column transformations are applied consistently to the objective, right-hand sides, and variable bounds, and an additional scalar balances the objective magnitude.  All subsequent formulas use the scaled data and omit the bars.  The primal and outer dual variables are mapped back to the original coordinates when the algorithm terminates.

\subsection{Outer updates}
\label{subsec:practical-outer}

The practical outer loop uses row-wise penalties instead of one scalar penalty.  The equality penalties are collected in $\sigma_{\E}^k\in\R^{m_{\E}}$, and the inequality penalties are collected in $\sigma_{\I}^k\in\R^{m_{\I}}$; both vectors have positive entries.  After each inner solve, the method updates the primal and outer dual variables, checks the original QP, and forms the next inner QP.  

The row-wise penalty update, proximal-center update, and infeasibility tests follow the corresponding mechanisms in QPALM \cite{hermans2022qpalm}.  Because \PDNQP{} solves each subproblem by a first-order method, it uses slower penalty growth and more gradual reduction of the inner stopping tolerances.  Anderson acceleration is added to improve the final stage of outer convergence.  \Cref{alg:outer} summarizes the complete outer iteration.

\begin{algorithm}[H]
\caption{Practical \PDNQP{} outer iteration}
\label{alg:outer}
\begin{algorithmic}[1]
\Require Scaled QP \eqref{eq:qp}; scalar proximal step $\gamma>0$; vectors $x^0,z^0\in\R^n$; vectors $\lambda_{\E}^0\in\R^{m_{\E}}$ and $\lambda_{\I}^0\in\R^{m_{\I}}$; row-penalty vectors $\sigma_{\E}^0\in\R^{m_{\E}}$ and $\sigma_{\I}^0\in\R^{m_{\I}}$ with positive entries; scalar tolerances $\varepsilon_a,\varepsilon_r,\delta_{a,0},\delta_{r,0},\varepsilon_{a,0},\varepsilon_{r,0}$; outer-iteration limit $K_{\mathrm{out}}$; parameters in \cref{tab:parameters}
\For{$k=0,1,\ldots,K_{\mathrm{out}}-1$}
    \State Form the inner QP \eqref{eq:inner-qp} using $z^k$, $\lambda_{\E}^k$, $\lambda_{\I}^k$, $\sigma_{\E}^k$, and $\sigma_{\I}^k$
    \State Apply \cref{alg:inner} and obtain an inner primal vector $w^{k+1}$ and an inner-solver status
    \If{the inner-solver status is numerical error}
        \State \Return numerical error
    \EndIf
    \If{the inner termination criteria are not satisfied}
        \State Retain the current outer variables and continue the inner solve
        \State \textbf{continue}
    \EndIf
    \State Extract the primal variable $x^{k+1}$ and update $s_{\I}^{k+1}$, $\lambda_{\E}^{k+1}$, and $\lambda_{\I}^{k+1}$ by \eqref{eq:practical-multiplier-update}
    \State Evaluate the outer residuals \eqref{eq:outer-residuals} and the termination criteria \eqref{eq:outer-termination}
    \If{the requested outer accuracy is reached}
        \State Unscale and \Return $(x^{k+1},\lambda_{\E}^{k+1},\lambda_{\I}^{k+1})$
    \EndIf
    \State Check the infeasibility conditions in \cref{subsec:infeasibility}
    \If{an infeasibility condition is accepted}
        \State \Return the corresponding status
    \EndIf
    \State Update the row penalties according to \cref{subsec:penalty-update}
    \State Update the inner stopping tolerances according to \cref{subsec:inner-tolerance-update}
    \State Update the proximal center and its tolerances according to \cref{subsec:center-update}
    \State Apply the Anderson step in \cref{subsec:anderson} when its safeguards hold
\EndFor
\State \Return iteration limit
\end{algorithmic}
\end{algorithm}

\subsubsection{Update of the primal and dual variables ($x^{k+1},s_{\I}^{k+1},\lambda_{\E}^{k+1},\lambda_{\I}^{k+1}$)}
\label{subsec:outer-residuals}

To control the progress of each coordinate, the method maintains positive row-penalty vectors $\sigma_{\E}^k$ and $\sigma_{\I}^k$ for the equality and inequality constraints, respectively.
The primal variable $x^{k+1}$ is obtained from the inner QP. The slack variable
and the dual variables are updated by
\begin{subequations}
\label{eq:practical-multiplier-update}
\begin{align}
 s_{\I}^{k+1}
 &=\max\!\left\{
 b_{\I},
 A_{\I}x^{k+1}+\lambda_{\I}^k\oslash\sigma_{\I}^k
 \right\},\\
 \lambda_{\E}^{k+1}
 &=\lambda_{\E}^k
   +\sigma_{\E}^k\odot(A_{\E}x^{k+1}-b_{\E}),\\
 \lambda_{\I}^{k+1}
 &=\min\!\left\{
 0,
 \lambda_{\I}^k
 +\sigma_{\I}^k\odot(A_{\I}x^{k+1}-b_{\I})
 \right\}\ ,
\end{align}
\end{subequations}
where the symbols $\odot$, $\oslash$, $\max$, $\min$ denote componentwise multiplication, division, maximization and minimization, respectively.

\subsubsection{Termination check}
\label{subsec:outer-termination}
The termination criteria follow from QPALM \cite{hermans2022qpalm}. In particular, we measure primal infeasibility and projected stationarity, namely, 
\begin{subequations}
\label{eq:outer-residuals}
\begin{align}
 r_{\mathrm p}^{k+1}
 &=\max\!\left\{
 \norminf{A_{\E}x^{k+1}-b_{\E}},
 \norminf{(b_{\I}-A_{\I}x^{k+1})_+}
 \right\},\\
 r_{\mathrm d}^{k+1}
 &=\norminf{
 x^{k+1}-\proj_{[\ell,u]}\!\left(
 x^{k+1}-Qx^{k+1}-c
 -A_{\E}^{\top}\lambda_{\E}^{k+1}
 -A_{\I}^{\top}\lambda_{\I}^{k+1}
 \right)}.
\end{align}
\end{subequations}
The primal residual checks the equality and violated inequality rows, while the projected residual accounts for stationarity under the variable bounds.

Let $s_{\mathrm p}^{k+1}$ and $s_{\mathrm d}^{k+1}$ denote the scales of the constraint and stationarity terms, respectively:
\begin{subequations}
\label{eq:outer-scales}
\begin{align}
 s_{\mathrm p}^{k+1}
 &=\max\!\left\{
 1,\norminf{A_{\E}x^{k+1}},\norminf{b_{\E}},
 \norminf{A_{\I}x^{k+1}},\norminf{b_{\I}}
 \right\},\\
 s_{\mathrm d}^{k+1}
 &=\max\!\left\{
 1,\norminf{Qx^{k+1}},\norminf{c},
 \norminf{A_{\E}^{\top}\lambda_{\E}^{k+1}
          +A_{\I}^{\top}\lambda_{\I}^{k+1}}
 \right\}.
\end{align}
\end{subequations}
Given the user-provided absolute tolerance $\varepsilon_a$ and relative tolerance $\varepsilon_r$, the algorithm terminates when
\begin{equation}
\label{eq:outer-termination}
 r_{\mathrm p}^{k+1}
 \le\varepsilon_a+\varepsilon_r s_{\mathrm p}^{k+1},
 \qquad
 r_{\mathrm d}^{k+1}
 \le\varepsilon_a+\varepsilon_r s_{\mathrm d}^{k+1}.
\end{equation}

\subsubsection{Row-wise penalty update ($\sigma_{\E}^k,\sigma_{\I}^k$)}
\label{subsec:penalty-update}

The penalty of a constraint row is increased when its primal residual does not decrease enough, while the penalties of the other rows are unchanged.

For a row $i$, let $r_i^k$ denote its primal residual at outer iteration $k$.  For $i\in\E$, $r_i^k=(A_{\E}x^k-b_{\E})_i$, and for $i\in\I$, $r_i^k=(A_{\I}x^k-s_{\I}^k)_i$.  The penalty of an equality row is increased using the same rule as in QPALM \cite{hermans2022qpalm}, when
\begin{equation}
\label{eq:penalty-trigger-eq}
 |r_i^{k+1}|>0.25|r_i^k|.
\end{equation}
The same test is used for an inequality row when the row is active.

Compared with QPALM, \PDNQP{} increases the penalties more slowly because the inner problem is solved by a first-order method.  A large penalty can make the inner QP more ill-conditioned and slow down the inner solver.

\subsubsection{Update of the inner stopping tolerances ($\delta_{a,k},\delta_{r,k}$)}
\label{subsec:inner-tolerance-update}

The tolerances are loose in early outer iterations and are reduced as the outer iteration approaches stationarity.  In particular, the reduction factor is chosen by
\begin{equation}
\label{eq:inner-tolerance-factor}
\rho_k^{\mathrm{in}}
=
\rho_{\max}^{\mathrm{in}}
-
\left(\rho_{\max}^{\mathrm{in}}-\rho_{\min}^{\mathrm{in}}\right)
\min\!\left\{
1,\,
\frac{\gamma^{-1}\norminf{x^{k+1}-z^k}}
     {r_{\mathrm d}^{k+1}}
\right\}.
\end{equation}
Here, $\rho_{\min}^{\mathrm{in}}=0.3$ and
$\rho_{\max}^{\mathrm{in}}=0.95$. $\rho_k^{\mathrm{in}}$ takes a smaller value when the proximal term is large relative to $r_{\mathrm d}^{k+1}$, and a larger value otherwise. The inner tolerances are then multiplied by $\rho_k^{\mathrm{in}}$.

Compared with QPALM, \PDNQP{} uses values of $\rho_k^{\mathrm{in}}$ closer to one.  In early outer iterations, when the proximal residual is small relative to the projected-stationarity residual, solving the inner QP to high accuracy gives little benefit because the outer progress is still dominated by constraint infeasibility. As the proximal residual becomes more important, the inner solve is tightened. This avoids spending many first-order iterations on unnecessarily accurate inner solves.

\subsubsection{Update of the proximal center and its tolerances ($z^{k+1},\varepsilon_{a,k},\varepsilon_{r,k}$)}
\label{subsec:center-update}
The center is updated to the new primal point only when the current feasibility threshold is satisfied.  Using the outer primal residual and scale in \eqref{eq:outer-residuals} and \eqref{eq:outer-scales}, the update is
\begin{equation}
\label{eq:center-update}
 z^{k+1}=
 \begin{cases}
 x^{k+1},
 &\text{if }r_{\mathrm p}^{k+1}
 \le\varepsilon_{a,k}+\varepsilon_{r,k}s_{\mathrm p}^{k+1},\\
 z^k,&\text{otherwise}.
 \end{cases}
\end{equation}

After the center is updated, the center-update tolerances are reduced by
\begin{equation}
\label{eq:center-tolerance-update}
 \varepsilon_{a,k+1}=\rho_z\varepsilon_{a,k},
 \qquad
 \varepsilon_{r,k+1}=\rho_z\varepsilon_{r,k},
\end{equation}
where $\rho_z=0.95$ is the reduction factor.  This lets the method move the proximal center with loose feasibility in early iterations, while requiring more accurate feasibility later.

\subsubsection{Infeasibility detection}
\label{subsec:infeasibility}
For primal infeasibility, \PDNQP{} checks whether the change in the outer dual variables gives an approximate Farkas certificate.  The candidate direction must satisfy the required sign and stationarity conditions and must separate the constraint system from zero.  If these conditions hold, the original QP is declared primal infeasible.

For dual infeasibility, the method checks whether the change in the primal variables gives an approximate recession direction.  The direction must satisfy the equality constraints, the inequality directions, and the finite variable bounds.  The problem is declared dual infeasible if the objective decreases along this direction through either negative quadratic curvature or a negative linear term.

\subsubsection{Anderson acceleration of the outer dual variables and proximal center ($\lambda_{\E},\lambda_{\I},z$)}
\label{subsec:anderson}

Unlike QPALM, \PDNQP{} applies Anderson acceleration to the outer dual variables and proximal center after primal feasibility and the row penalties become stable. At this stage, the outer iteration usually enters a more stable local convergence regime, while the remaining progress can still be slow. Anderson acceleration is therefore applied to improve the tail convergence.

Let $m_A=5$ be the Anderson memory size.  The method stores the most recent $m_A$ outer iterates of $(\lambda_{\E},\lambda_{\I},z)$.  Classical Type-II Anderson acceleration \cite{anderson1965iterative} is then applied to these iterates to form an extrapolated point for the next outer iteration.

\subsection{Inner solver}
\label{subsec:practical-inner}

Each outer iteration forms a strongly convex QP.  In \eqref{eq:inner-qp}, these QPs share the same sparse constraint matrix $B$, while their penalty terms and linear terms can change between outer iterations.

The practical inner solver follows the rAPDHG method in \cref{subsec:base-pdhg} and \cite{lu2025practicalqp}.  Unlike QPALM, the inner solver is designed for a first-order method and therefore uses a different numerical treatment.  We first warm start and scale the inner QP, then set the primal weight and diagonal step sizes. The inner iterations use adaptive step acceptance, restarting, and the inner termination check described below. \Cref{alg:inner} summarizes the inner solver.

\begin{algorithm}[H]
\caption{Practical diagonal PDHG method for the inner QP}
\label{alg:inner}
\begin{algorithmic}[1]
\Require Convex inner QP \eqref{eq:inner-qp}; initial vectors $w^0\in\R^{d_w}$ and $y^0\in\R^{m_{\E}+m_{\I}}$; scalar primal weight $\omega>0$; scalar inner stopping tolerances $\delta_{a,k},\delta_{r,k}$; accepted-iteration limit $K_{\mathrm{in}}$
\State Warm start the inner variables according to \cref{subsec:warm-start}
\State Apply the inner scaling in \cref{subsec:inner-scaling} and form the diagonal primal and dual steps by \eqref{eq:diagonal-scales}--\eqref{eq:diagonal-steps}
\For{$t=0,1,\ldots,K_{\mathrm{in}}-1$}
    \Repeat
        \State Compute the trial inner primal--dual point $(w^{\mathrm{tr}},y^{\mathrm{tr}})$ by \eqref{eq:pdqp-trial}
        \If{the adaptive condition \eqref{eq:adaptive-condition} fails}
            \State Reduce the diagonal primal and dual steps and reject the trial
        \Else
            \State Accept $(w^{t+1},y^{t+1})\gets(w^{\mathrm{tr}},y^{\mathrm{tr}})$ and update the weighted averages
        \EndIf
    \Until{the trial is accepted or the rejection limit is reached}
    \If{the rejection limit is reached}
        \State \Return numerical error
    \EndIf
    \If{a stopping check is due}
        \State Test the current and averaged primal candidates according to \cref{subsec:inner-stopping}
        \If{either candidate satisfies the inner termination criteria}
            \State \Return the accepted inner primal--dual candidate
        \EndIf
    \EndIf
    \If{a restart condition in \cref{subsec:inner-restart} holds}
        \State Select the better current-or-average candidate and update the primal weight according to \cref{subsec:primal-weight}
        \State Rebuild the diagonal steps and clear the weighted averages
    \EndIf
\EndFor
\State \Return the current inner primal--dual variables with iteration-limit status
\end{algorithmic}
\end{algorithm}

\subsubsection{Warm start of the inner variables}
\label{subsec:warm-start}

Successive outer iterations form closely related inner QPs.  \PDNQP{} therefore uses the solution of the previous inner QP to initialize the next one.  When the row penalties change, the residual and slack variables are updated to match the new inner QP.

Using the current primal variable $x$, inequality dual variable $\lambda_{\I}$, and inequality penalty vector $\sigma_{\I}$, the residual and slack variables are initialized by
\begin{equation}
\label{eq:inner-warm-start}
 s_{\I}^{0}=\max\!\left\{
 b_{\I},A_{\I}x+\lambda_{\I}\oslash\sigma_{\I}
 \right\},
 \qquad
 r_{\E}^{0}=A_{\E}x-b_{\E},
 \qquad
 r_{\I}^{0}=A_{\I}x-s_{\I}^{0}.
\end{equation}
The second and third expressions restore the residual definitions in \eqref{eq:prox-alm-subproblem}, while the first updates the inequality split using the new row penalties.  This keeps the previous primal progress while making the initial inner variables consistent with the new subproblem.

\subsubsection{Scaling of the inner primal and dual variables ($w,y$)}
\label{subsec:inner-scaling}

The row-wise penalties can have very different values across different outer iterations.  Since $\sigma_{\E,i}^k$ and $\sigma_{\I,i}^k$ appear in the diagonal residual blocks of $H_k$, large differences in these penalties can make the inner QP poorly scaled.  \PDNQP{} first scales each residual variable by its own row penalty:
\begin{equation}
\label{eq:inner-penalty-scaling}
 r_{\E,i}\leftarrow \sqrt{\sigma_{\E,i}^k}\,r_{\E,i},
 \qquad
 r_{\I,i}\leftarrow \sqrt{\sigma_{\I,i}^k}\,r_{\I,i}.
\end{equation}
This removes the scale difference in the residual blocks caused by the row-wise penalties.

After this penalty scaling, \PDNQP{} applies $8$ iterations of Ruiz scaling \cite{ruiz2001scaling} to the inner QP.  As in the outer scaling, positive diagonal row and column scaling matrices are used to balance the variables and constraints.  The same transformations are applied to the Hessian, linear term, right-hand side, and variable bounds.  All inner iterations use the scaled problem.

\subsubsection{Inner termination check ($\delta_{a,k},\delta_{r,k}$)}
\label{subsec:inner-stopping}

The inner termination check uses the same form as the outer termination check in \eqref{eq:outer-residuals}--\eqref{eq:outer-termination}.  For an inner primal--dual candidate $(w,y)$, define
\begin{subequations}
\label{eq:inner-residuals}
\begin{align}
 r_{\mathrm p}^{\mathrm{in}}
 &=\norminf{Bw-h},\\
 r_{\mathrm d}^{\mathrm{in}}
 &=\norminf{
 w-\proj_{\mathcal W}\!\left(
 w-H_kw-q_k+B^\top y
 \right)}.
\end{align}
\end{subequations}
Here, $H_k$ and $q_k$ are the Hessian and linear term of the inner QP, and $\mathcal W=[\underline w,\overline w]$ is its box constraint.  Thus $r_{\mathrm p}^{\mathrm{in}}$ measures inner primal feasibility and $r_{\mathrm d}^{\mathrm{in}}$ measures projected stationarity.

The corresponding primal and dual scales are
\begin{subequations}
\label{eq:inner-scales}
\begin{align}
 s_{\mathrm p}^{\mathrm{in}}
 &=\max\!\left\{
 1,\norminf{Bw},\norminf{h}
 \right\},\\
 s_{\mathrm d}^{\mathrm{in}}
 &=\max\!\left\{
 1,\norminf{H_kw},\norminf{q_k},
 \norminf{B^\top y}
 \right\}.
\end{align}
\end{subequations}
The inner solver returns when
\begin{equation}
\label{eq:inner-termination}
 r_{\mathrm p}^{\mathrm{in}}
 \le \delta_{a,k}+\delta_{r,k}s_{\mathrm p}^{\mathrm{in}},
 \qquad
 r_{\mathrm d}^{\mathrm{in}}
 \le \delta_{a,k}+\delta_{r,k}s_{\mathrm d}^{\mathrm{in}}.
\end{equation}

\subsubsection{Primal-weight update ($\omega$)}
\label{subsec:primal-weight}

The primal weight $\omega>0$ balances the primal and dual steps in rAPDHG.  Following PDLP \cite{applegate2021pdlp}, it is updated at each restart using the primal and dual movements from the previous restart.

Let $j$ denote the restart epoch.  Define
\begin{equation}
\label{eq:primal-weight-movement}
 \Delta_w^j=\normtwo{w^{j,0}-w^{j-1,0}},
 \qquad
 \Delta_y^j=\normtwo{y^{j,0}-y^{j-1,0}}.
\end{equation}
When both changes are nonzero, \PDNQP{} updates the primal weight by
\begin{equation}
\label{eq:primal-weight-update}
 \omega^j
 =
 \exp\!\left(
 0.25\log\frac{\Delta_y^j}{\Delta_w^j}
 +0.75\log\omega^{j-1}
 \right).
\end{equation}
The ratio $\Delta_y^j/\Delta_w^j$ estimates the primal weight that balances the primal and dual movements, while the log-scale averaging avoids large changes between restarts.

\subsubsection{Diagonal step sizes for the inner primal and dual variables ($\tau,\eta$)}
\label{subsec:diagonal-preconditioner}

The base rAPDHG method uses one primal step size and one dual step size.  For the practical inner solver, different rows and columns can have very different scales.  \PDNQP{} therefore uses a separate primal step $\tau_i>0$ for each component $w_i$ and a separate dual step $\eta_j>0$ for each constraint row.

Let $H$ and $B$ denote the scaled inner Hessian and constraint matrix.  For the primal component $w_i$, let $H_{:i}$ and $B_{:i}$ be the $i$th columns of $H$ and $B$.  Let $B_{ji}$ be entry $(j,i)$ of $B$, and let $\sigma_j>0$ be the row penalty corresponding to constraint row $j$.  We define
\begin{subequations}
\label{eq:diagonal-scales}
\begin{align}
 m_{w,i}
 &=
 \norm{H_{:i}}_1+
 \begin{cases}
 \displaystyle
 \max\!\left\{
 \norm{B_{:i}}_1,
 \left(\sum_j\sigma_jB_{ji}^2\right)^{1/2}
 \right\},
 &w_i\text{ belongs to }x\text{ or }s_{\I},\\[2mm]
 \norm{B_{:i}}_1,
 &w_i\text{ belongs to }r_{\E}\text{ or }r_{\I},
 \end{cases}\\
 m_{y,j}
 &=\sum_{i:m_{w,i}>0}\frac{B_{ji}^2}{m_{w,i}}.
\end{align}
\end{subequations}
Here, $\norm{H_{:i}}_1$ measures the size of the quadratic term for $w_i$, and $\norm{B_{:i}}_1$ measures its coupling with the constraints.  For the $x$ and $s_{\I}$ blocks, the penalty-weighted term also accounts for rows with large penalties.  The quantity $m_{y,j}$ then scales dual row $j$ using the primal scales of the variables in that row.

Let $\omega>0$ be the primal weight used in rAPDHG.  The diagonal step sizes are
\begin{equation}
\label{eq:diagonal-steps}
 \tau_i=\frac{0.99}{\omega m_{w,i}},
 \qquad
 \eta_j=\frac{0.99\,\omega}{m_{y,j}}.
\end{equation}
Thus a primal component or dual row with a larger scale receives a smaller step size.  The factor $0.99$ gives a small margin for the PDHG step.

\subsubsection{Adaptive quadratic PDHG step ($w,y$)}
\label{subsec:pdqp-step}

The diagonal step sizes give an initial step for each primal and dual component.  \PDNQP{} further checks each trial step before accepting it.  The trial update follows the accelerated PDHG iteration in \cite{lu2025practicalqp}, using the diagonal step vectors $\tau$ and $\eta$.

Let $w^t$ and $\bar w^t$ be the current and averaged inner primal variables, and let $y^t$ be the current inner dual variable.  Let $N_p^t$ and $N_d^t$ be the numbers of accepted primal and dual steps, and set
$\alpha_t=2/(N_p^t+2)$ and $\beta_t=N_d^t/(N_d^t+1)$.  The trial point is
\begin{subequations}
\label{eq:pdqp-trial}
\begin{align}
 w^{\mathrm{tr}}
 &=\proj_{\mathcal W}\!\left(
 w^t-\tau\odot
 \left[
 q+(1-\alpha_t)H\bar w^t
 +\alpha_tHw^t-B^\top y^t
 \right]
 \right),\\
 y^{\mathrm{tr}}
 &=y^t+\eta\odot
 \left[
 h-Bw^{\mathrm{tr}}
 -\beta_t(Bw^{\mathrm{tr}}-Bw^t)
 \right].
\end{align}
\end{subequations}
The primal step uses the accelerated quadratic gradient and then projects onto $\mathcal W$.  The dual step uses the new constraint residual and the change from the previous primal point.

Set $\Delta w=w^{\mathrm{tr}}-w^t$ and $\Delta y=y^{\mathrm{tr}}-y^t$.  The trial is accepted when
\begin{equation}
\label{eq:adaptive-condition}
 \frac12\sum_i\frac{(\Delta w_i)^2}{\tau_i}
 +\frac12\sum_j\frac{(\Delta y_j)^2}{\eta_j}
 \ge
 \frac12\sum_i\left|\Delta w_i(H\Delta w)_i\right|
 +\sum_i\left|\Delta w_i(B^\top\Delta y)_i\right|.
\end{equation}
The left-hand side measures the primal and dual changes relative to their step sizes.  The two terms on the right-hand side measure the quadratic change from $H$ and the primal--dual coupling through $B$.  The trial step uses this condition to check whether the current diagonal step sizes are suitable for the local QP before updating the iterate.  This allows the inner solver to adapt its step sizes to the local quadratic and constraint structure.

\subsubsection{Restart}
\label{subsec:inner-restart}

The restart rule follows PDLP and PDQP
\cite{applegate2021pdlp,lu2025practicalqp}.  At each restart, the method compares the current and averaged primal--dual iterates and uses the better one as the starting point of the next epoch.  This removes old iterates from the average and lets the new epoch follow the current progress.

\subsection{GPU implementation}
\label{subsec:gpu}

The inner solver mainly uses sparse matrix--vector products, vector operations, projections, and reductions.  \PDNQP{} keeps these repeated operations on the GPU and reuses the fixed sparse structure across outer iterations.

\paragraph{Fixed sparse matrices.}
The sparse structures of $H_k$, $B$, and their transposes are built once and stored on the GPU.  A row-penalty update changes only the diagonal residual blocks of $H_k$, while an update of the outer dual variables or proximal center changes only the linear term $q_k$.  Therefore, the sparse matrix structures are reused throughout the solve.

\paragraph{GPU inner iterations.}
The primal and dual updates, step acceptance, averaging, restart, and inner termination checks are performed on the GPU.  These repeated operations are captured by CUDA graphs to reduce CPU--GPU synchronization during the inner solve.

\paragraph{CPU outer iteration.}
The CPU controls the outer updates, including the row penalties, termination and infeasibility checks, and Anderson acceleration.  The required vector operations and reductions are still performed on the GPU.  Only the results needed by the outer logic are returned to the CPU.

\paragraph{Algorithmic parameters.}
The default implementation parameters are listed in \cref{tab:parameters}.

\section{Numerical experiments}
\label{sec:experiments}

This section evaluates \PDNQP{} along three dimensions: (i) robustness on established nonconvex QP benchmarks, (ii) sensitivity to the requested solution accuracy, and (iii) scalability to million-variable instances.  \Cref{subsec:experimental-protocol} describes the test sets, solver configurations, stopping criteria, and computational environment.  \Cref{subsec:main-results} presents the computational results.  The construction of the million-variable instances and the complete instance-level results are provided in the appendix.

\subsection{Datasets and experimental setup}
\label{subsec:experimental-protocol}

\paragraph{Datasets.}
We consider two collections of nonconvex QPs.  The first consists of the $107$ CUTEst instances used in the QPALM study \cite{gould2015cutest,hermans2022qpalm}.  This collection provides an established benchmark for assessing robustness on nonconvex QPs of moderate size.
The second collection consists of $30$ synthetic, application-based QPs with approximately $10^6$ variables each.  It contains $10$ graph soft-labeling instances, $10$ concave-cost multistage-flow instances, and $10$ indefinite kernel-SVM instances.  These problems are designed to test scalability in regimes where sparse matrix factorizations can become prohibitively expensive.  The underlying models and instance-generation procedures are described in \cref{app:million-models,app:million-generation}.

\paragraph{Compared solvers.}
We compare \PDNQP{} with four established solvers for nonconvex QP or general nonlinear programming:

\begin{itemize}\item QPALM 1.2.6 \cite{hermans2022qpalm};\item Ipopt 3.14.19, accessed through cyipopt 1.7.0 \cite{wachter2006ipopt};\item MadNLP 0.9.2 with MadNLPGPU 0.8.0, CUDA.jl 5.11.3, and CUDSS.jl 0.6.5 \cite{shin2024gpu}; and\item Knitro 16.0.0 \cite{byrd2006knitro}.\end{itemize}

QPALM uses its nonconvex proximal mode and factorization method $2$, with $\rho=0.1$, $\sigma_{\mathrm{init}}=20$, $\delta=100$, $10$ scaling sweeps, and $\gamma_{\mathrm{init}}=\gamma_{\max}=10^7$.  Its initial inner tolerances are $1$, and its infeasibility tolerances are $\min\{\varepsilon_a,10^{-7}\}$.  Its formulation includes one identity constraint row for each variable.
Ipopt uses sequential MUMPS 5.8.2.  Its overall, feasibility, dual, and complementarity tolerances are set to the current solver tolerance.  MadNLPGPU uses a bound-relaxation factor of $10^{-8}$ and solves a sparse condensed KKT system with cuDSS.  Knitro uses its Interior/Direct algorithm with the exact quadratic Hessian, one local solve, and the supplied starting point.  Its feasibility and optimality tolerances are set to the current solver tolerance.  \PDNQP{} uses the default parameters stated in Appendix \ref{app:parameters}.

\paragraph{Initialization.} All solvers start from $x^0=0$ on the CUTEst instances. For each synthetic instance, all solvers use the same fixed starting point specified by the instance-generation procedure. \Cref{tab:benchmark-profiles} summarizes the three benchmark profiles.

\begin{table}[H]
\centering
\caption{The three benchmark profiles.  Every instance has a $3600$-second time limit.}
\label{tab:benchmark-profiles}
\begin{tabular}{@{}llrrr@{}}
\toprule
Profile & Test set and start & Cases & $\varepsilon_a=\varepsilon_r$ & Limit (s) \\
\midrule
CUTEst standard & CUTEst, $x^0=0$ & 107 & $10^{-4}$ & 3600 \\
CUTEst high accuracy & CUTEst, $x^0=0$ & 107 & $10^{-6}$ & 3600 \\
Million-variable & Frozen $x^0$ & 30 & $10^{-4}$ & 3600 \\
\bottomrule
\end{tabular}
\end{table}

\paragraph{Computing environment.}
The GPU experiments are performed on an NVIDIA H100 with 80~GB of HBM3 memory, and the CPU experiments use an Intel Xeon Platinum 8481C processor.  \PDNQP{} and MadNLPGPU run on the GPU, whereas QPALM, Ipopt, and Knitro run on the CPU.  All solvers use double-precision arithmetic. 

\paragraph{Termination and tolerance calibration.}

We use the termination criteria defined in \cref{subsec:outer-termination}. Because solvers use different internal stopping criteria, identical input tolerances do not necessarily yield comparable solution accuracy. We therefore evaluate every finite returned point using a common external KKT test applied to the original, unscaled formulation \eqref{eq:qp}. This test uses the primal feasibility and projected-stationarity residuals, together with the absolute-plus-relative scaling rule defined in \cref{subsec:outer-termination}. The primal feasibility test also includes variable-bound violations.

For each instance, the initial solver tolerance is set to the tolerance specified by the benchmark profile. If the returned point passes the external test, it is accepted. Otherwise, the solver tolerance is halved and the problem is solved again. The external acceptance tolerances remain fixed throughout this calibration process. The process stops at the first accepted point. An outcome is recorded as a failure (F) if an attempt reaches the $3600$-second time limit or if the returned point still fails the external test when the solver tolerance reaches $1/128$ of its initial value.

\paragraph{Timing and aggregate metric.}
We report solver-call times and assign a time of $3600$ seconds to every failure. For a benchmark profile with $N$ instances and assigned times $t_1,\ldots,t_N$, we report the shifted geometric mean
\begin{equation}
\label{eq:sgm10}
\operatorname{SGM10}(t_1,\ldots,t_N)
=\exp\left(\frac{1}{N}\sum_{i=1}^{N}\log(t_i+10)\right)-10.
\end{equation}
The $10$-second shift reduces the influence of very short runs. Lower SGM10 values indicate better aggregate performance.

\subsection{Results}
\label{subsec:main-results}

\Cref{tab:aggregate-results} summarizes the success counts, failure counts, and SGM10 values for the three benchmark profiles. Complete instance-level results are provided in Appendix \ref{app:instance-results}.

\begin{table}[ht]
\centering
\caption{Aggregate results.  S counts points accepted by the external KKT test and PI or DI classifications.  F counts all other outcomes.  Every F is assigned $3600$ seconds in SGM10.}
\label{tab:aggregate-results}
\small
\begin{tabular}{@{}llrrr@{}}
\toprule
Profile & Solver & S & F & $\operatorname{SGM10}$ (s) \\
\midrule
CUTEst $10^{-4}$ (107)
& \textbf{\PDNQP{}} & \textbf{107} & \textbf{0} & \textbf{4.4} \\
& QPALM       & 101 & 6  & 23.0 \\
& Ipopt       & 95  & 12 & 26.3 \\
& MadNLPGPU   & 101 & 6  & 11.3 \\
& Knitro      & 89  & 18 & 26.4 \\
\addlinespace
CUTEst $10^{-6}$ (107)
& \textbf{\PDNQP{}} & \textbf{101} & \textbf{6} & 44.3 \\
& QPALM       & 95 & 12 & \textbf{37.7} \\
& Ipopt       & 77 & 30 & 85.1 \\
& MadNLPGPU   & 86 & 21 & 37.9 \\
& Knitro      & 73 & 34 & 81.0 \\
\addlinespace
Million $10^{-4}$ (30)
& \textbf{\PDNQP{}} & \textbf{27} & \textbf{3} & \textbf{76.0} \\
& QPALM       & 0 & 30 & 3600.0 \\
& Ipopt       & 0 & 30 & 3600.0 \\
& MadNLPGPU   & 0 & 30 & 3600.0 \\
& Knitro      & 0 & 30 & 3600.0 \\
\bottomrule
\end{tabular}
\end{table}

\paragraph{Low accuracy.}
At a tolerance of $10^{-4}$, \PDNQP{} returns $104$ points that pass the external KKT test and reports two primal-infeasible classifications and one dual-infeasible classification. It is the only solver to achieve a successful outcome on all $107$ CUTEst instances. QPALM and MadNLPGPU each succeed on $101$ instances, followed by Ipopt on $95$ and Knitro on $89$. \PDNQP{} also achieves the lowest SGM10 of $4.4$ seconds, compared with $11.3$ seconds for MadNLPGPU, $23.0$ for QPALM, $26.3$ for Ipopt, and $26.4$ for Knitro.

\paragraph{High accuracy.}
At a tolerance of $10^{-6}$, \PDNQP{} returns $98$ points that pass the external KKT test and reports the same three infeasibility classifications. It retains the highest success count, with $101$ successful outcomes, followed by QPALM with $95$, MadNLPGPU with $86$, Ipopt with $77$, and Knitro with $73$. QPALM and MadNLPGPU achieve lower SGM10 values of $37.7$ and $37.9$ seconds, respectively, compared with $44.3$ seconds for \PDNQP{}, but succeed on $6$ and $15$ fewer instances. Thus, \PDNQP{} maintains the highest success rate at the tighter tolerance, although it no longer achieves the lowest SGM10.

\paragraph{Effect of solution accuracy.}
Tightening the tolerance from $10^{-4}$ to $10^{-6}$ increases the $\operatorname{SGM10}$ of \PDNQP{} from $4.4$ to $44.3$ seconds, by about one order of magnitude.  This behavior is expected for a first-order method. At high accuracy, the inexpensive inner iterations must be repeated many more times to reduce the final KKT residuals.  In contrast, Newton-type and interior-point methods can benefit from faster local convergence once they are close to a solution. Therefore, \PDNQP{} has a clearer advantage at lower accuracy, while its first-order convergence becomes more costly at $10^{-6}$.

\paragraph{Million-variable instances.}
\PDNQP{} passes the external KKT test on $27$ of the $30$ million-variable instances and has an SGM10 of $76.0$ seconds.  QPALM, Ipopt, MadNLPGPU, and Knitro do not solve any instance in this profile.  \PDNQP{} is therefore the only tested solver that returns verified solutions at this scale.  These results demonstrate the scalability of the factorization-free \PDNQP{} framework on million-variable nonconvex QPs.

\paragraph{Dependence on problem structure.}
\PDNQP{} performs best on problems such as NCVXQP, QPNBAND, STNQP, and the flow and SVM instances. These problems have simple sparse coupling: each variable interacts with only a small number of other variables or constraint rows, and many constraints are simple bounds. In this case, the inner PDHG iterations remain cheap and the residual form adds little extra coupling. \PDNQP{} is less competitive on the MPC and soft-labeling instances, where many general constraints couple the variables and the residual form introduces more auxiliary variables. These results suggest that \PDNQP{} is particularly suitable for QPs with sparse local quadratic and constraint structure, but is less favorable when the problem contains many strongly coupled general constraints.

\section{Conclusion}

We presented \PDNQP{}, a factorization-free first-order method for large-scale nonconvex quadratic programming. The method combines a proximal augmented-Lagrangian framework with a residual-form reformulation that preserves a fixed sparse constraint operator across the strongly convex inner QPs. These subproblems are solved by restarted accelerated PDHG using sparse matrix--vector products, projections, and vector operations, which makes the method suitable for GPU computation.

The numerical results show that \PDNQP{} is robust on the nonconvex CUTEst benchmark and is particularly effective at lower accuracy. Its advantage is strongest when the quadratic and constraint structures have sparse local coupling, so that the inner first-order iterations remain inexpensive. On the million-variable test set, \PDNQP{} returns verified solutions on $27$ of $30$ instances, while none of the other tested solvers succeeds. These results show that factorization-free first-order methods can provide a practical approach for nonconvex QPs whose structure becomes difficult for sparse factorization.

\bibliographystyle{amsplain}
\bibliography{ref-papers}

\appendix

\section{Default algorithmic parameters}
\label{app:parameters}
\label{subsec:defaults}

The defaults used by the practical solver are listed below.

\begin{longtable}{@{}p{0.28\textwidth}p{0.16\textwidth}p{0.48\textwidth}@{}}
\caption{Default parameters of \PDNQP{}.}
\label{tab:parameters}\\
\toprule
Parameter & Default & Role \\
\midrule
\endfirsthead
\toprule
Parameter & Default & Role \\
\midrule
\endhead

$\varepsilon_a,\varepsilon_r$
& $10^{-4},10^{-4}$
& Outer absolute and relative KKT tolerances. \\

$K_{\mathrm{out}}$
& $2^{31}-1$
& Outer iteration limit. \\

Time limit
& $+\infty$
& Wall-clock termination limit. \\

KKT pass limit
& $+\infty$
& Cumulative inner matrix-pass limit. \\

$K_{\mathrm{in}}$
& $40000$
& Accepted inner-iteration limit per outer call. \\

$\delta_{a,0},\delta_{r,0},\varepsilon_{a,0},\varepsilon_{r,0}$
& $1,1,1,1$
& Initial inner stopping and center-update tolerances. \\

$\sigma_{\mathrm{init}}$
& $20$
& Initial row-penalty scale. \\

$\theta_{\sigma}$
& $0.25$
& Residual reduction threshold for the row-penalty update. \\

$\delta_{\gamma}$
& $10^{-2}$
& Margin used in the proximal regularization. \\

$\rho_{\min}^{\mathrm{in}},\rho_{\max}^{\mathrm{in}}$
& $0.3,0.95$
& Fast and slow reduction factors for the inner tolerances. \\

$\rho_z$
& $0.95$
& Reduction factor for the center-update tolerances. \\

$m_A$
& $5$
& Anderson memory size. \\

Ruiz scaling sweeps
& $8$
& Number of outer and inner scaling sweeps. \\

\bottomrule
\end{longtable}

\section{Million-variable benchmark construction}
\label{app:million-models}

The 30 application-based QPs contain ten instances from each of three families.  Their requested size is $10^6$ variables; minor rounding is needed when variables are indexed by classes or network layers.  For every family, scenarios $s=101,\ldots,110$ were generated once, validated, and frozen before the experiments.

\paragraph{Soft-label graph clustering.}
Let $G=(V,E)$ be a graph with positive edge weights $w_{ij}$.  Let $K$ be the number of classes.  The full decision vector is $x\in\R^{|V|K}$.  For vertex $i$, the block $x_i=(x_{i1},\ldots,x_{iK})\in\R^K$ contains its class memberships.  The scalar $a_{ic}$ is the unary cost of assigning vertex $i$ to class $c$.  The scalar $\delta\in(0,1)$ controls the allowed class sizes.  Define the weighted degree $d_i=\sum_j w_{ij}$ and the normalized edge weight $\widetilde w_{ij}=w_{ij}/\sqrt{d_i d_j}$.  The model is
\begin{align*}
 \min_x\quad&
 -\sum_{\{i,j\}\in E}\widetilde w_{ij}\sum_{c=1}^{K}x_{ic}x_{jc}
 +\sum_{i\in V}\sum_{c=1}^{K}a_{ic}x_{ic},
 \\
 \text{s.t.}\quad&
 \sum_{c=1}^{K}x_{ic}=1 \quad (i\in V),\qquad 0\le x_{ic}\le1,\nonumber\\
 &\frac{(1-\delta)\lvert V\rvert}{K}
 \le\sum_{i\in V}x_{ic}
 \le\frac{(1+\delta)\lvert V\rvert}{K}
 \quad (c=1,\ldots,K).\nonumber
\end{align*}
Thus $x_i$ is a probability vector.  Graph affinity rewards neighboring vertices that share a label, while the class-capacity interval prevents collapse to a single class.  A balanced planted labeling generates the graph and unary costs.  With $N=\lfloor10^6/K\rfloor$, the instance has $n=NK$, $N$ simplex equalities, and $2K$ class-capacity inequalities.  The frozen start is $x_{ic}^0=1/K$, which satisfies every constraint.

\paragraph{Concave-cost multistage flow.}
Let the network have $n_v$ nodes and $n_e$ arcs.  Let $B\in\R^{n_v\times n_e}$ be the node--arc incidence matrix.  Let $f\in\R^{n_e}$ be the flow vector and let $b\in\R^{n_v}$ be the supply vector.  For arc $e$, the scalars $u_e>0$ and $c_e>0$ are its capacity and linear unit cost.  The scalar $\alpha\in(0,\frac{1}{2})$ controls the concave discount.  The model is
\begin{align*}
 \min_f\quad&
 \sum_e\left(c_e f_e-\alpha\frac{c_e}{u_e}f_e^2\right),
 \\
 \text{s.t.}\quad& Bf=b,\qquad 0\le f_e\le u_e.\nonumber
\end{align*}
The Hessian is diagonal with $Q_{ee}=-2\alpha c_e/u_e<0$.  The condition $\alpha<1/2$ keeps the marginal cost positive over the capacity interval, so the model represents increasing cost with economies of scale.  The network has $L$ layers of width $w$ and three carriers.  Seeded permutations connect adjacent layers, giving $n=3w(L-1)$ arcs.  A reference flow $f^0$ defines $b=Bf^0$, so the frozen start is exactly feasible.

\paragraph{Indefinite graph-kernel SVM.}
Let $N$ be the number of samples.  The label vector is $y\in\{-1,1\}^N$, with balanced labels.  Let $Y=\diag(y)\in\R^{N\times N}$.  The decision vector is $\alpha\in\R^N$, and $U_\alpha>0$ is its componentwise upper bound.  Let $\mathbf 1\in\R^N$ be the all-ones vector.  Let $K_\eta\in\R^{N\times N}$ be the graph-kernel matrix.  The model is
\begin{align*}
 \min_{\alpha}\quad&
 \tfrac{1}{2}\alpha^\top YK_\eta Y\alpha-\mathbf{1}^\top\alpha,
 \\
 \text{s.t.}\quad&
 y^\top\alpha=0,\qquad 0\le\alpha_i\le U_\alpha.\nonumber
\end{align*}
Each sample has a feature vector $z_i\in\R^d$.  Let $\xi_i\in\R^d$ satisfy $\xi_i\sim\mathcal N(0,I_d)$, and let $\zeta_i\in\R^{d-1}$ satisfy $\zeta_i\sim\mathcal N(0,I_{d-1})$.  Let $e_1\in\R^d$ be the first coordinate vector.  The scalars $\nu>0$ and $\chi>0$ are the noise scale and class separation.  The feature generator is $z_i=\nu\xi_i+\chi y_i e_1+0.15(0,\zeta_i)$.  Let $k$ be the number of nearest neighbors and let $\kappa>0$ be the kernel scale.  Let $\eta\in(1/2,1)$ be the kernel mixing parameter.  Let $W\in\R^{N\times N}$ be the graph weight matrix, and set $D=\diag(W\mathbf 1)\in\R^{N\times N}$.  Let $S,K_\eta\in\R^{N\times N}$ be the normalized graph and kernel matrices.  A symmetric $k$-nearest-neighbor graph defines
\[
 W_{ij}=\left\lvert\tanh\!\left(\frac{\kappa z_i^\top z_j}{d}\right)\right\rvert,
 \qquad
 S=D^{-\frac{1}{2}}WD^{-\frac{1}{2}},\qquad
 K_\eta=(1-\eta)I+\eta S.
\]
Sufficiently negative eigenvalues of $S$ make $K_\eta$ indefinite.  Congruence by $Y$ preserves its inertia.  Each instance has $10^6$ variables and one equality.  The frozen start is $\alpha^0=0\in\R^N$.

\section{Frozen million-variable scenario distributions}
\label{app:million-generation}

For reproducibility, let $a,b\in\R$ with $a<b$.  The notation \(\mathcal U(a,b)\) denotes a continuous uniform draw.  All finite sets are sampled uniformly.  A separate NumPy \texttt{default\_rng} stream is initialized by \texttt{SeedSequence([s,r])}, with family code \(r=1,2,3\) for soft-label, flow, and SVM, respectively.

\begin{table}[H]
\centering
\small
\caption{Scenario distributions for the frozen million-variable instances.}
\label{tab:million-generation}
\begin{tabularx}{\textwidth}{@{}lX@{}}
\toprule
Family & Parameters drawn for each scenario seed \(s\) \\
\midrule
Soft-label &
\(K\in\{3,5,8,10\}\), sampled degree \(r_g\in\{8,\ldots,24\}\),
same-label neighbor fraction \(\mathcal U(0.55,0.90)\), cross-class base
similarity \(\mathcal U(0.08,0.45)\), log-weight standard deviation
\(\mathcal U(0.10,0.50)\), anchor fraction \(\mathcal U(0.03,0.20)\),
anchor strength \(\mathcal U(0.5,2)\), unary-noise standard deviation
\(\mathcal U(0.01,0.08)\), and \(\delta\sim\mathcal U(0.05,0.25)\). \\
\addlinespace
Flow &
\(L\in\{4,6,8,12,16\}\), \(\alpha\sim\mathcal U(0.03,0.24)\),
base cost \(\mathcal U(0.05,0.30)\), capacity ratio
\(\rho\sim\mathcal U(1.05,1.80)\), and stage-growth rate
\(\mathcal U(0,0.12)\).  The geometric carrier-cost ratio is drawn between
\((1-2\alpha)^{-1}+0.25\) and \((1-2\alpha)^{-1}+2\).
Bounds for the carrier reference flows are drawn from
\(\mathcal U(0.30,1.00)\) and \(\mathcal U(1.20,3.00)\).  For each carrier,
\(h_k\) is then drawn uniformly within that interval. \\
\addlinespace
SVM &
\(d\in\{4,8,12,16\}\), \(k\in\{6,8,12,16,24\}\), class separation
\(\mathcal U(0.75,3.00)\), feature-noise standard deviation
\(\mathcal U(0.50,1.25)\), kernel scale \(\kappa\sim\mathcal U(0.75,3.00)\),
\(\eta\sim\mathcal U(0.60,0.97)\), and
\(U_\alpha\sim\mathcal U(0.25,2.00)\). \\
\bottomrule
\end{tabularx}
\end{table}

\section{Per-instance results}
\label{app:instance-results}

\input{pdnqp_instance_results.tex}

\end{document}

%% file: pdnqp_instance_results.tex
Tables~\ref{tab:cutest-standard-detail}--\ref{tab:million-detail} report every instance.  The CUTEst columns $n$ and $m$ are the dimensions stored in the frozen metadata.  For the million-variable set, $m=m_{\E}+m_{\I}$.  Variable bounds are excluded.  Time is solver-call time, and Obj is evaluated in the original formulation.  For a point accepted after tolerance reduction, Time is the solver-call time of the first accepted attempt; earlier attempts are not included.  A superscript PI or DI denotes primal or dual infeasibility, and the preceding number is the measured solver-call time.  F denotes every other outcome and is assigned $3600$ seconds when computing $\operatorname{SGM10}$.

\begingroup
\scriptsize
\setlength{\tabcolsep}{1.0pt}
\begin{longtable}{@{}lcc|cc|cc|cc|cc|cc@{}}
\caption{Per-instance CUTEst results for $\varepsilon_a=\varepsilon_r=10^{-4}$.  The common start is $x^0=0$, and each solver attempt has a $3600$-second limit.}\label{tab:cutest-standard-detail}\\
\hline
Problem & $n$ & $m$ & \multicolumn{2}{|c}{\PDNQP{}} & \multicolumn{2}{|c}{QPALM} & \multicolumn{2}{|c}{Ipopt} & \multicolumn{2}{|c}{MadNLPGPU} & \multicolumn{2}{|c}{Knitro} \\
\cline{4-13}
& & & Time & Obj & Time & Obj & Time & Obj & Time & Obj & Time & Obj \\
\hline
\endfirsthead
\multicolumn{13}{c}{\tablename\ \thetable\ (continued)} \\
\hline
Problem & $n$ & $m$ & \multicolumn{2}{|c}{\PDNQP{}} & \multicolumn{2}{|c}{QPALM} & \multicolumn{2}{|c}{Ipopt} & \multicolumn{2}{|c}{MadNLPGPU} & \multicolumn{2}{|c}{Knitro} \\
\cline{4-13}
& & & Time & Obj & Time & Obj & Time & Obj & Time & Obj & Time & Obj \\
\hline
\endhead
\hline\multicolumn{13}{r}{Continued on next page}\\
\endfoot
\hline
\endlastfoot
A0ENDNDL & 45006 & 15002 & 1.8e+00 & 0.0e+00 & 3.8e+00 & 1.3e-02 & 5.7e-01 & 4.5e-02 & 1.2e+00 & 1.1e-02 & 5.2e-01 & 1.3e-05 \\
A0ENINDL & 45006 & 15002 & 1.7e+00 & 0.0e+00 & 5.0e+00 & -2.9e-02 & 5.7e-01 & 4.5e-02 & 5.0e+01 & 4.5e-02 & 6.1e-01 & 1.5e-05 \\
A0ENSNDL & 45006 & 15002 & 1.1e+00 & 0.0e+00 & 4.5e+00 & -1.3e-02 & 3.6e+01 & 3.6e-02 & 2.6e+00 & 6.1e-01 & 6.7e-01 & 8.4e-06 \\
A0ESDNDL & 45006 & 15002 & 2.4e+00 & 0.0e+00 & 6.4e+00 & -1.7e-02 & 5.5e-01 & 4.5e-02 & 1.2e+00 & 4.5e-02 & 4.5e-01 & 1.4e-05 \\
A0ESINDL & 45006 & 15002 & 2.3e+00 & 0.0e+00 & 5.2e+00 & 8.6e-03 & 8.5e-01 & 4.5e-02 & 1.3e+00 & 4.5e-02 & 4.9e-01 & 1.6e-05 \\
A0ESSNDL & 45006 & 15002 & 6.8e-01 & 0.0e+00 & 3.1e+00 & -4.7e-05 & 1.8e+01 & 3.6e-02 & 2.4e+00 & 6.1e-01 & 6.9e-01 & 8.0e-06 \\
A0NNDNDL & 60012 & 20004 & 3.3e+01 & 0.0e+00 & 1.7e+02 & -1.3e-03 & 1.7e+00 & 4.5e-02 & 5.8e+00 & 4.5e-02 & 8.8e-01 & 1.2e-05 \\
A0NNDNIL & 60012 & 20004 & 7.6e+01 & 0.0e+00 & 8.3e+02 & 8.3e-02 & 1.2e+01 & 4.5e-02 & 8.1e+01 & 7.5e-01 & 6.9e+00 & 4.8e-02 \\
A0NNDNSL & 60012 & 20004 & 2.6e+00 & 0.0e+00 & 1.5e+01 & 2.4e-03 & 5.4e+00 & 4.1e-02 & 6.6e+00 & 1.3e+01 & 1.1e+00 & 2.1e-03 \\
A0NNSNSL & 60012 & 20004 & 2.0e+00 & 0.0e+00 & 1.1e+01 & 2.3e-05 & 8.9e+00 & 3.6e-02 & 1.1e+01 & 6.0e-01 & 9.5e-01 & 9.1e-03 \\
A0NSDSDL & 60012 & 20004 & 1.2e+01 & 0.0e+00 & 2.3e+01 & 9.0e-05 & 1.8e+00 & 4.5e-02 & 7.5e+00 & 2.2e-02 & 4.1e-01 & 7.3e-07 \\
A0NSDSDS & 6012 & 2004 & 6.8e-01 & 0.0e+00 & 8.7e-01 & -7.8e-06 & 7.4e-01 & 3.4e-03 & 7.5e-01 & 1.5e+00 & 2.6e-01 & 5.6e-06 \\
A0NSDSIL & 60012 & 20004 & 2.1e+01 & 0.0e+00 & 1.1e+02 & -1.9e-02 & 1.2e+01 & 4.5e-02 & 1.8e+02 & 7.5e-01 & 4.4e+00 & 2.8e-05 \\
A0NSDSSL & 60012 & 20004 & 2.9e+00 & 0.0e+00 & 1.6e+01 & -8.8e-02 & 2.8e+00 & 3.7e-02 & 8.8e+00 & 9.3e+01 & 9.6e-01 & 2.7e-05 \\
A0NSSSSL & 60012 & 20004 & 5.7e+00 & 0.0e+00 & 1.7e+01 & -1.5e-02 & 5.8e+00 & 3.4e-02 & 1.1e+01 & 8.1e+01 & 7.7e-01 & 1.3e-03 \\
A2ENDNDL & 45006 & 15002 & 1.2e+00 & 4.9e-05 & 7.7e+00 & 3.2e-05 & 1.1e+00 & 1.9e-03 & 1.9e+00 & 2.0e-03 & 5.1e-01 & 1.7e-04 \\
A2ENINDL & 45006 & 15002 & 7.7e-01 & 4.6e-05 & 9.1e+00 & 1.8e-05 & 1.4e+00 & 1.9e-03 & 1.9e+00 & 1.4e-03 & 5.3e-01 & 1.8e-04 \\
A2ENSNDL & 45006 & 15002 & 7.8e-01 & 1.3e-04 & 7.1e+00 & -4.3e-05 & 2.0e+01 & 3.6e-02 & 2.5e+00 & 2.7e+00 & 6.4e-01 & 2.4e-02 \\
A2ESDNDL & 45006 & 15002 & 1.5e+00 & 7.0e-05 & 9.3e+00 & 5.2e-06 & 1.2e+00 & 1.9e-03 & 1.9e+00 & 2.0e-03 & 5.2e-01 & 1.8e-04 \\
A2ESINDL & 45006 & 15002 & 8.8e-01 & 3.7e-05 & 9.3e+00 & 1.6e-05 & 9.4e-01 & 1.9e-03 & 1.8e+00 & 1.4e-03 & 5.3e-01 & 1.8e-04 \\
A2ESSNDL & 45006 & 15002 & 8.9e-01 & 1.1e-04 & 4.7e+00 & 1.9e-04 & 2.4e+01 & 4.2e-02 & 2.6e+00 & 2.7e+00 & 6.4e-01 & 2.4e-02 \\
A2NNDNDL & 60012 & 20004 & 1.5e+00 & 2.5e-04 & 2.2e+01 & 1.0e-03 & 2.7e+00 & 1.5e-03 & 6.9e+00 & 5.5e-03 & 1.1e+00 & 2.6e-04 \\
A2NNDNIL & 60012 & 20004 & 7.0e+00 & 5.5e+01 & 1.9e+01 & 6.0e-01 & F & -- & F & -- & 4.5e+00 & 1.5e+00 \\
A2NNDNSL & 60012 & 20004 & 2.2e+00 & 3.4e-04 & 1.6e+01 & 3.4e-05 & 5.6e+00 & 2.0e-02 & 5.8e+00 & 1.3e+01 & 1.0e+00 & 2.0e-03 \\
A2NNSNSL & 60012 & 20004 & 3.4e+00 & 0.0e+00 & 1.6e+01 & 5.5e-05 & 2.1e+01 & -5.8e-02 & 9.9e+00 & 7.7e+01 & 8.6e-01 & 1.2e-02 \\
A2NSDSDL & 60012 & 20004 & 4.3e+00 & 3.3e-06 & 3.1e+01 & 5.1e-05 & 2.1e+00 & 1.0e-03 & 9.2e+00 & 8.5e-04 & 7.1e-01 & 7.4e-04 \\
A2NSDSIL & 60012 & 20004 & 7.9e+00 & 1.6e+01 & 2.8e+01 & 8.7e+00 & 2.2e+01 & 1.9e+00 & 1.3e+02 & 1.0e+00 & 5.6e+00 & 5.8e+00 \\
A2NSDSSL & 60012 & 20004 & 4.2e+00 & 3.2e-08 & 2.6e+01 & 5.4e-03 & 5.0e+00 & -5.1e-01 & 8.8e+00 & 9.1e+01 & 8.4e-01 & 1.5e-04 \\
A2NSSSSL & 60012 & 20004 & 4.2e+00 & 0.0e+00 & 1.7e+01 & 6.4e-03 & 1.0e+01 & 3.5e-02 & 1.0e+01 & 3.8e+02 & 8.2e-01 & 3.3e-04 \\
A5ENDNDL & 45006 & 15002 & 6.7e-01 & 1.2e-04 & 9.4e+00 & 3.3e-05 & 9.6e-01 & 3.0e-03 & 1.9e+00 & 3.4e-03 & 8.2e-01 & 9.5e-04 \\
A5ENINDL & 45006 & 15002 & 8.3e-01 & 1.6e-05 & 1.7e+01 & 4.0e-05 & 1.0e+00 & 3.0e-03 & 2.1e+00 & 2.9e-03 & 5.3e-01 & 1.0e-03 \\
A5ENSNDL & 45006 & 15002 & 9.2e-01 & 2.6e-05 & 5.3e+00 & -3.9e-05 & 1.6e+01 & 8.2e-02 & 2.3e+00 & 1.5e+01 & 7.2e-01 & 1.0e-01 \\
A5ESDNDL & 45006 & 15002 & 1.4e+00 & 1.7e-04 & 7.6e+00 & 8.2e-05 & 9.3e-01 & 3.0e-03 & 2.0e+00 & 3.4e-03 & 5.3e-01 & 9.7e-04 \\
A5ESINDL & 45006 & 15002 & 9.1e-01 & 2.3e-06 & 9.9e+00 & 1.2e-04 & 1.2e+00 & 3.0e-03 & 1.8e+00 & 2.9e-03 & 5.2e-01 & 1.1e-03 \\
A5ESSNDL & 45006 & 15002 & 9.4e-01 & 2.9e-05 & 5.9e+00 & -3.5e-03 & 1.5e+01 & 3.9e-02 & 2.4e+00 & 1.5e+01 & 6.3e-01 & 1.0e-01 \\
A5NNDNDL & 60012 & 20004 & 1.5e+00 & 2.8e-05 & 4.5e+01 & 1.3e-03 & 3.3e+00 & 2.7e-03 & 8.7e+00 & 5.4e-03 & 1.1e+00 & 3.0e-04 \\
A5NNDNIL & 60012 & 20004 & 8.1e+00 & 4.1e+01 & 2.4e+01 & 1.0e+01 & F & -- & F & -- & 4.6e+00 & 3.8e+00 \\
A5NNDNSL & 60012 & 20004 & 2.3e+00 & 1.8e-04 & 1.8e+01 & 3.4e-05 & 5.0e+00 & -5.3e-01 & 1.5e+01 & 3.7e-02 & 1.4e+00 & 2.4e-03 \\
A5NNSNSL & 60012 & 20004 & 5.5e+00 & 0.0e+00 & 1.9e+01 & -1.8e-05 & 3.5e+01 & -8.2e-01 & 9.1e+00 & 7.8e+01 & 7.5e-01 & 5.4e-03 \\
A5NSDSDL & 60012 & 20004 & 4.5e+00 & 5.2e-06 & 5.5e+01 & 1.4e-04 & 2.0e+00 & 2.7e-03 & 1.0e+01 & 1.6e-03 & 8.2e-01 & 4.3e-04 \\
A5NSDSDM & 6012 & 2004 & 7.1e-01 & 0.0e+00 & 8.9e-01 & -7.8e-06 & 5.5e-01 & 3.4e-03 & 8.0e-01 & 1.5e+00 & 3.1e-01 & 5.6e-06 \\
A5NSDSIL & 60012 & 20004 & 5.8e+00 & 1.9e+01 & 2.3e+01 & 1.4e+01 & 3.1e+01 & 1.3e+00 & 1.0e+02 & 1.6e+00 & 3.5e+00 & 1.8e+00 \\
A5NSDSSL & 60012 & 20004 & 4.1e+00 & 1.1e-06 & 2.5e+01 & -6.5e-04 & 8.9e+00 & 3.9e-02 & 9.4e+00 & 1.3e+01 & 8.6e-01 & 9.1e-03 \\
A5NSSNSM & 6012 & 2004 & 7.0e-01 & 0.0e+00 & 7.4e-01 & -7.8e-06 & 5.5e-01 & 3.4e-03 & 7.6e-01 & 1.5e+00 & 3.9e-01 & 5.6e-06 \\
A5NSSSSL & 60012 & 20004 & 3.5e+00 & 1.6e-09 & 2.2e+01 & 7.5e-03 & 1.3e+01 & 3.9e-02 & 1.1e+01 & 1.4e+01 & 8.5e-01 & 6.6e-03 \\
BIGGSC4 & 4 & 13 & 5.4e-02 & -2.4e+01 & 7.6e-04 & -2.4e+01 & 1.6e-02 & -2.4e+01 & 3.2e+00 & -2.4e+01 & 1.8e-01 & -2.4e+01 \\
BLOCKQP1 & 10010 & 5001 & 1.4e+00 & 5.0e+00 & 1.1e-01 & -5.0e+03 & 6.3e-01 & 5.0e+00 & F & -- & 3.2e-01 & 5.0e+00 \\
BLOCKQP2 & 10010 & 5001 & 7.3e-01 & -5.0e+03 & 1.9e-01 & -5.0e+03 & 9.6e-01 & -5.0e+03 & 1.7e+00 & -5.0e+03 & 5.7e-01 & -5.0e+03 \\
BLOCKQP3 & 10010 & 5001 & 5.4e+00 & 5.0e+00 & 1.1e-01 & 5.0e+00 & F & -- & F & -- & 3.9e-01 & 5.0e+00 \\
BLOCKQP4 & 10010 & 5001 & 6.6e-01 & -2.5e+03 & 3.3e-01 & -2.5e+03 & 9.0e-01 & -2.5e+03 & 1.3e+00 & -2.5e+03 & 6.6e-01 & -2.5e+03 \\
BLOCKQP5 & 10010 & 5001 & 1.5e+00 & 5.0e+00 & 1.1e-01 & 5.0e+00 & F & -- & F & -- & 3.9e-01 & 5.0e+00 \\
BLOWEYA & 4002 & 2002 & 8.3e-01 & -1.7e-04 & 1.3e-02 & -1.2e-03 & 4.3e+01 & -1.6e-02 & 1.3e-01 & -2.2e-02 & 2.5e-01 & -2.1e-02 \\
BLOWEYB & 4002 & 2002 & 4.3e-01 & -3.4e-05 & 1.7e-02 & -7.3e-04 & 2.3e+01 & -3.8e-03 & 7.9e-02 & -5.6e-06 & 2.5e-01 & 0.0e+00 \\
BLOWEYC & 22 & 12 & 2.1e-01 & -3.1e+00 & 1.4e-03 & -3.1e+00 & 4.9e-03 & -3.1e+00 & 2.2e-02 & -3.1e+00 & 1.8e-01 & -3.1e+00 \\
CLEUVEN3 & 1200 & 2973 & 8.5e-01 & 4.0e+06 & 1.6e+00 & 1.6e+07 & 2.0e+01 & 2.9e+05 & 1.3e+00 & 2.9e+05 & 2.3e+00 & 2.9e+05 \\
CLEUVEN4 & 1200 & 2973 & 3.1e+00 & 4.1e+06 & 1.5e+02 & 1.2e+07 & 5.7e+01 & 2.9e+05 & 3.7e+01 & 2.9e+05 & 2.3e+01 & 2.9e+05 \\
CLEUVEN5 & 1200 & 2973 & 8.6e-01 & 4.0e+06 & 1.6e+00 & 1.6e+07 & 2.1e+01 & 2.9e+05 & 1.2e+00 & 2.9e+05 & 2.4e+00 & 2.9e+05 \\
CLEUVEN6 & 1200 & 3091 & 4.2e+00 & 7.0e+07 & 1.7e+00 & 2.3e+07 & 6.0e+01 & 2.2e+07 & 5.7e-01 & 2.2e+07 & 7.5e+00 & 2.2e+07 \\
FERRISDC & 2200 & 210 & 8.7e-02 & -1.4e-12 & 5.9e-03 & 0.0e+00 & 3.0e+00 & -1.5e-05 & 2.6e-01 & -2.8e-05 & 2.3e-01 & 0.0e+00 \\
GOULDQP1 & 32 & 17 & 9.2e-01 & -3.5e+03 & 2.1e-03 & -3.5e+03 & 1.7e-02 & -3.5e+03 & 1.0e-01 & -3.5e+03 & 1.9e-01 & -3.5e+03 \\
HATFLDH & 4 & 13 & 5.1e-02 & -2.4e+01 & 5.6e-04 & -2.4e+01 & 8.4e-03 & -2.4e+01 & 7.4e-02 & -2.4e+01 & 2.1e-01 & -2.4e+01 \\
HS44 & 4 & 6 & 9.5e-02 & -1.3e+01 & 5.3e-04 & -1.5e+01 & 1.0e-02 & -1.3e+01 & 5.7e-02 & -1.3e+01 & 2.1e-01 & -1.5e+01 \\
HS44NEW & 4 & 6 & 9.5e-02 & -1.3e+01 & 4.8e-04 & -1.5e+01 & 8.7e-03 & -1.3e+01 & 8.7e-02 & -1.3e+01 & 1.8e-01 & -1.5e+01 \\
LEUVEN2 & 1530 & 2329 & 3.3e+01 & -1.4e+07 & 1.8e+00 & -1.4e+07 & 2.6e+00 & -1.4e+07 & 6.8e-01 & -1.4e+07 & F & -- \\
LEUVEN3 & 1200 & 2973 & 1.2e+00 & -5.2e+08 & 1.2e+03 & -1.7e+09 & 3.8e+02 & -2.0e+09 & 7.3e+02 & -2.0e+09 & 7.7e+01 & -1.9e+09 \\
LEUVEN4 & 1200 & 2973 & 3.2e+00 & -5.0e+08 & F & -- & 1.6e+03 & -2.0e+09 & 2.1e+01 & -2.0e+09 & 4.7e+02 & -2.0e+09 \\
LEUVEN5 & 1200 & 2973 & 1.2e+00 & -5.2e+08 & 1.2e+03 & -1.7e+09 & 3.8e+02 & -2.0e+09 & 5.7e+00 & -2.0e+09 & 7.7e+01 & -1.9e+09 \\
LEUVEN6 & 1200 & 3091 & 4.8e+00 & -1.8e+08 & 1.4e+03 & -1.2e+09 & 1.6e+02 & -1.2e+09 & 3.9e+00 & -1.2e+09 & 1.1e+02 & -9.6e+08 \\
LEUVEN7 & 360 & 946 & 1.5e+01 & 6.9e+02 & 5.4e-02 & 6.9e+02 & 3.1e-01 & 6.9e+02 & 1.0e-01 & 6.9e+02 & 3.4e-01 & 6.9e+02 \\
LINCONT & 1257 & 419 & 1.4e-01$^{\mathrm{PI}}$ & -- & 1.1e-01$^{\mathrm{PI}}$ & -- & 1.1e+00$^{\mathrm{PI}}$ & -- & 9.1e-01$^{\mathrm{PI}}$ & -- & 5.5e+00$^{\mathrm{PI}}$ & -- \\
MPC1 & 2550 & 3833 & 1.7e+01 & -2.3e+07 & 2.4e+00 & -2.3e+07 & 6.0e+00 & -2.3e+07 & 1.3e+00 & -2.3e+07 & F & -- \\
MPC2 & 1530 & 2351 & 1.1e+01 & -1.5e+07 & 8.3e-01 & -1.5e+07 & 7.5e-01 & -1.5e+07 & 5.1e-01 & -1.5e+07 & F & -- \\
MPC3 & 1530 & 2351 & 1.5e+01 & -1.5e+07 & 1.0e+00 & -1.5e+07 & 8.6e-01 & -1.5e+07 & 5.1e-01 & -1.5e+07 & F & -- \\
MPC4 & 1530 & 2351 & 1.3e+01 & -1.5e+07 & 8.8e-01 & -1.5e+07 & 6.6e-01 & -1.5e+07 & 3.4e-01 & -1.5e+07 & F & -- \\
MPC5 & 1530 & 2351 & 1.1e+01 & -1.5e+07 & 1.0e+00 & -1.5e+07 & 7.7e-01 & -1.5e+07 & 3.9e-01 & -1.5e+07 & F & -- \\
MPC6 & 1530 & 2351 & 1.5e+01 & -1.5e+07 & 1.1e+00 & -1.5e+07 & 8.3e-01 & -1.5e+07 & 3.1e-01 & -1.5e+07 & F & -- \\
MPC7 & 1530 & 2351 & 1.5e+01 & -1.5e+07 & 1.0e+00 & -1.5e+07 & 8.5e-01 & -1.5e+07 & 4.3e-01 & -1.5e+07 & F & -- \\
MPC8 & 1530 & 2351 & 1.3e+01 & -1.5e+07 & 1.1e+00 & -1.5e+07 & 9.0e-01 & -1.5e+07 & 3.3e-01 & -1.5e+07 & F & -- \\
MPC9 & 1530 & 2351 & 1.7e+01 & -1.5e+07 & 9.8e-01 & -1.5e+07 & 1.1e+00 & -1.5e+07 & 3.2e-01 & -1.5e+07 & F & -- \\
MPC10 & 1530 & 2351 & 1.7e+01 & -1.5e+07 & 1.1e+00 & -1.5e+07 & 1.0e+00 & -1.5e+07 & 2.5e-01 & -1.5e+07 & F & -- \\
MPC11 & 1530 & 2351 & 1.7e+01 & -1.5e+07 & 1.1e+00 & -1.5e+07 & 8.3e-01 & -1.5e+07 & 5.8e-01 & -1.5e+07 & F & -- \\
MPC12 & 1530 & 2351 & 1.4e+01 & -1.5e+07 & 1.2e+00 & -1.5e+07 & 7.5e-01 & -1.5e+07 & 6.6e-01 & -1.5e+07 & F & -- \\
MPC13 & 1530 & 2351 & 1.9e+01 & -1.5e+07 & 1.0e+00 & -1.5e+07 & 7.8e-01 & -1.5e+07 & 3.8e-01 & -1.5e+07 & F & -- \\
MPC14 & 1530 & 2351 & 1.2e+01 & -1.5e+07 & 1.1e+00 & -1.5e+07 & 8.2e-01 & -1.5e+07 & 3.7e-01 & -1.5e+07 & F & -- \\
MPC15 & 1530 & 2351 & 1.3e+01 & -1.5e+07 & 1.1e+00 & -1.5e+07 & 8.2e-01 & -1.5e+07 & 3.7e-01 & -1.5e+07 & F & -- \\
MPC16 & 1530 & 2351 & 1.3e+01 & -1.5e+07 & 1.0e+00 & -1.5e+07 & 8.8e-01 & -1.5e+07 & 4.8e-01 & -1.5e+07 & F & -- \\
NASH & 72 & 24 & 1.5e-02$^{\mathrm{PI}}$ & -- & 1.0e-03$^{\mathrm{PI}}$ & -- & 2.1e-02$^{\mathrm{PI}}$ & -- & 1.4e-01$^{\mathrm{PI}}$ & -- & 0.0e+00$^{\mathrm{PI}}$ & -- \\
NCVXQP1 & 10000 & 5000 & 5.4e+00 & -7.5e+09 & 2.5e+03 & -7.5e+09 & F & -- & 1.3e+01 & -7.5e+09 & 1.9e+01 & -7.5e+09 \\
NCVXQP2 & 10000 & 5000 & 3.6e+00 & -5.8e+09 & F & -- & F & -- & 1.9e+01 & -5.8e+09 & 6.6e+01 & -5.8e+09 \\
NCVXQP3 & 10000 & 5000 & 5.8e+00 & -3.1e+09 & F & -- & F & -- & 2.4e+01 & -3.0e+09 & 5.7e+01 & -3.1e+09 \\
NCVXQP4 & 10000 & 2500 & 4.6e-01 & -9.4e+09 & 8.8e+02 & -9.4e+09 & 1.5e+03 & -9.4e+09 & 1.2e+01 & -9.4e+09 & 3.1e+01 & -9.4e+09 \\
NCVXQP5 & 10000 & 2500 & 7.0e-01 & -6.6e+09 & 1.6e+03 & -6.6e+09 & 1.9e+03 & -6.6e+09 & 1.4e+01 & -6.6e+09 & 3.8e+01 & -6.6e+09 \\
NCVXQP6 & 10000 & 2500 & 8.6e-01 & -3.3e+09 & 2.2e+03 & -3.4e+09 & 2.7e+03 & -3.4e+09 & 3.2e+01 & -3.4e+09 & 2.8e+01 & -3.5e+09 \\
NCVXQP7 & 10000 & 7500 & 2.0e+01 & -5.2e+09 & F & -- & F & -- & 6.7e+00 & -5.2e+09 & 1.2e+01 & -5.2e+09 \\
NCVXQP8 & 10000 & 7500 & 1.4e+01 & -3.6e+09 & F & -- & F & -- & 1.4e+01 & -3.6e+09 & 2.4e+01 & -3.6e+09 \\
NCVXQP9 & 10000 & 7500 & 8.4e+00 & -2.1e+09 & 3.3e+03 & -2.1e+09 & F & -- & 1.4e+01 & -2.1e+09 & 6.0e+01 & -2.1e+09 \\
PORTSNQP & 10 & 2 & 3.8e-01 & 1.9e+00 & 5.2e-04 & 1.9e+00 & 4.8e-03 & 1.9e+00 & 3.0e-02 & 1.9e+00 & 1.7e-01 & 1.9e+00 \\
QPNBAND & 50000 & 25000 & 3.7e-01 & -2.5e+05 & 6.4e+00 & -2.5e+05 & 2.9e+00 & -2.5e+05 & 1.4e-01 & -2.5e+05 & 7.8e-01 & -2.5e+05 \\
QPNBLEND & 83 & 74 & 3.7e-01 & -9.2e-03 & 2.1e-03 & -9.1e-03 & 1.3e-02 & -9.1e-03 & 5.4e-02 & -9.1e-03 & 2.1e-01 & -9.1e-03 \\
QPNBOEI1 & 384 & 440 & 1.0e+00 & 6.8e+06 & 6.3e-02 & 6.8e+06 & 1.2e+00 & 6.8e+06 & 1.7e+00 & 6.8e+06 & 2.4e-01 & 6.8e+06 \\
QPNBOEI2 & 143 & 185 & 1.2e+00 & 1.4e+06 & 1.6e-02 & 1.4e+06 & 4.0e-01 & 1.4e+06 & 1.4e+00 & 1.4e+06 & 2.7e-01 & 1.4e+06 \\
QPNSTAIR & 467 & 356 & 8.9e-01 & 5.1e+06 & 3.8e-02 & 5.1e+06 & 6.7e-01 & 5.1e+06 & 8.2e-01 & 5.1e+06 & F & -- \\
SOSQP1 & 5000 & 2501 & 1.4e-01 & 2.4e-05 & 4.5e-02 & -1.5e-01 & 6.8e-02 & 7.3e-05 & 2.5e-01 & -2.0e-05 & 2.4e-01 & 0.0e+00 \\
SOSQP2 & 5000 & 2501 & 5.2e+00 & -1.2e+03 & F & -- & F & -- & F & -- & 2.7e-01 & -1.2e+03 \\
STATIC3 & 434 & 96 & 3.5e-02$^{\mathrm{DI}}$ & -- & 2.1e-03$^{\mathrm{DI}}$ & -- & F & -- & 2.3e-01$^{\mathrm{DI}}$ & -- & 1.9e-01$^{\mathrm{DI}}$ & -- \\
STNQP1 & 8193 & 4095 & 7.5e-01 & -3.1e+05 & 8.6e+00 & -3.1e+05 & 1.1e+02 & -3.1e+05 & 3.0e-01 & -3.1e+05 & 5.3e-01 & -3.1e+05 \\
STNQP2 & 8193 & 4095 & 6.2e-01 & -5.7e+05 & 2.7e+01 & -5.7e+05 & 5.7e+00 & -5.7e+05 & 4.1e-01 & -5.7e+05 & 6.4e-01 & -5.7e+05 \\
\end{longtable}
\endgroup

\begingroup
\scriptsize
\setlength{\tabcolsep}{1.0pt}
\begin{longtable}{@{}lcc|cc|cc|cc|cc|cc@{}}
\caption{Per-instance CUTEst results for $\varepsilon_a=\varepsilon_r=10^{-6}$.  The common start is $x^0=0$, and each solver attempt has a $3600$-second limit.}\label{tab:cutest-high-detail}\\
\hline
Problem & $n$ & $m$ & \multicolumn{2}{|c}{\PDNQP{}} & \multicolumn{2}{|c}{QPALM} & \multicolumn{2}{|c}{Ipopt} & \multicolumn{2}{|c}{MadNLPGPU} & \multicolumn{2}{|c}{Knitro} \\
\cline{4-13}
& & & Time & Obj & Time & Obj & Time & Obj & Time & Obj & Time & Obj \\
\hline
\endfirsthead
\multicolumn{13}{c}{\tablename\ \thetable\ (continued)} \\
\hline
Problem & $n$ & $m$ & \multicolumn{2}{|c}{\PDNQP{}} & \multicolumn{2}{|c}{QPALM} & \multicolumn{2}{|c}{Ipopt} & \multicolumn{2}{|c}{MadNLPGPU} & \multicolumn{2}{|c}{Knitro} \\
\cline{4-13}
& & & Time & Obj & Time & Obj & Time & Obj & Time & Obj & Time & Obj \\
\hline
\endhead
\hline\multicolumn{13}{r}{Continued on next page}\\
\endfoot
\hline
\endlastfoot
A0ENDNDL & 45006 & 15002 & 1.5e+00 & 0.0e+00 & 4.3e+00 & -6.4e-08 & 6.6e-01 & 1.8e-04 & 4.2e+01 & 1.8e-04 & 4.3e-01 & 1.3e-05 \\
A0ENINDL & 45006 & 15002 & 1.6e+00 & 0.0e+00 & 3.2e+00 & 1.0e-06 & 7.4e-01 & 1.8e-04 & 4.2e+01 & 1.8e-04 & 4.7e-01 & 1.5e-05 \\
A0ENSNDL & 45006 & 15002 & 1.3e+00 & 0.0e+00 & 3.4e+00 & -3.6e-06 & 1.2e+01 & 9.2e-05 & 2.2e+01 & 1.3e-04 & 5.9e-01 & 8.4e-06 \\
A0ESDNDL & 45006 & 15002 & 2.4e+00 & 0.0e+00 & 3.3e+00 & -3.5e-06 & 6.1e-01 & 1.8e-04 & 1.2e+00 & 1.8e-04 & 6.2e-01 & 1.4e-05 \\
A0ESINDL & 45006 & 15002 & 2.4e+00 & 0.0e+00 & 2.9e+00 & 7.3e-08 & 8.8e-01 & 1.8e-04 & 1.5e+00 & 1.8e-04 & 5.3e-01 & 1.6e-05 \\
A0ESSNDL & 45006 & 15002 & 8.4e-01 & 0.0e+00 & 3.3e+00 & 8.8e-07 & 1.8e+01 & 1.5e-04 & 4.0e+00 & 1.7e-04 & 5.8e-01 & 8.0e-06 \\
A0NNDNDL & 60012 & 20004 & 4.7e+02 & 0.0e+00 & 7.3e+01 & -8.2e-06 & 1.8e+00 & 1.8e-04 & 5.5e+00 & 1.8e-04 & 1.1e+00 & 1.2e-05 \\
A0NNDNIL & 60012 & 20004 & 1.5e+03 & 0.0e+00 & F & -- & 1.0e+01 & 2.0e-04 & 8.0e+01 & 1.9e-04 & 7.2e+00 & 4.4e-06 \\
A0NNDNSL & 60012 & 20004 & 1.5e+01 & 0.0e+00 & 1.4e+01 & 2.2e-04 & 6.7e+00 & -1.4e-03 & 1.3e+01 & 1.2e-04 & 6.2e-01 & 1.9e-07 \\
A0NNSNSL & 60012 & 20004 & 4.5e+00 & 0.0e+00 & 1.2e+01 & -3.5e-07 & 1.6e+01 & 2.0e-05 & 2.8e+01 & -1.0e-04 & 7.6e-01 & 6.1e-07 \\
A0NSDSDL & 60012 & 20004 & 6.2e+00 & 0.0e+00 & 2.0e+01 & -7.6e-06 & 1.5e+00 & 1.8e-04 & 7.2e+00 & 1.8e-04 & 7.9e-01 & 7.3e-07 \\
A0NSDSDS & 6012 & 2004 & 1.1e+00 & 0.0e+00 & 1.1e+00 & -2.4e-07 & 6.5e-01 & -2.6e-05 & 9.4e-01 & 1.9e-05 & 3.4e-01 & 1.1e-11 \\
A0NSDSIL & 60012 & 20004 & 9.1e+01 & 0.0e+00 & 8.0e+01 & 1.3e-04 & 1.3e+01 & 2.0e-04 & 2.4e+02 & 2.8e-04 & 4.5e+00 & 6.8e-08 \\
A0NSDSSL & 60012 & 20004 & 4.0e+00 & 0.0e+00 & 2.1e+01 & -1.7e-06 & 2.7e+00 & -3.9e-03 & 1.2e+01 & -7.0e-05 & 1.0e+00 & 5.4e-09 \\
A0NSSSSL & 60012 & 20004 & 3.8e+00 & 0.0e+00 & 2.2e+01 & 3.5e-05 & 8.9e+00 & -2.0e-03 & 1.5e+01 & -1.7e-03 & 9.4e-01 & 1.6e-07 \\
A2ENDNDL & 45006 & 15002 & 1.7e+00 & 3.2e-09 & 2.2e+01 & -5.9e-09 & F & -- & F & -- & F & -- \\
A2ENINDL & 45006 & 15002 & 4.4e+00 & 4.5e-09 & 2.6e+01 & -4.2e-08 & F & -- & F & -- & F & -- \\
A2ENSNDL & 45006 & 15002 & 1.3e+00 & 4.7e-10 & 5.5e+00 & 4.3e-09 & 1.7e+01 & 3.7e-04 & 4.1e+00 & 6.5e-04 & F & -- \\
A2ESDNDL & 45006 & 15002 & 2.2e+00 & 3.6e-09 & 2.4e+01 & -5.1e-08 & F & -- & F & -- & F & -- \\
A2ESINDL & 45006 & 15002 & 1.5e+00 & 2.1e-09 & 2.6e+01 & -1.4e-08 & F & -- & F & -- & F & -- \\
A2ESSNDL & 45006 & 15002 & 1.5e+00 & 3.2e-10 & 5.8e+00 & 4.3e-09 & 2.7e+01 & 2.5e-04 & 3.5e+00 & 3.3e-04 & 6.0e-01 & 3.6e-04 \\
A2NNDNDL & 60012 & 20004 & 3.0e+02 & 2.1e-06 & 2.1e+03 & 1.4e-06 & F & -- & F & -- & F & -- \\
A2NNDNIL & 60012 & 20004 & 1.2e+01 & 5.5e+01 & 3.6e+01 & 6.0e-01 & F & -- & F & -- & F & -- \\
A2NNDNSL & 60012 & 20004 & 5.1e+01 & 2.2e-08 & 1.5e+02 & 6.1e-08 & 1.0e+01 & -4.7e-03 & 7.9e+00 & 1.2e-04 & 9.0e-01 & 9.2e-07 \\
A2NNSNSL & 60012 & 20004 & 4.0e+00 & 6.2e-10 & 3.2e+01 & 7.0e-08 & 3.5e+02 & -3.0e-04 & 1.3e+01 & 1.1e-04 & 1.1e+00 & 6.2e-07 \\
A2NSDSDL & 60012 & 20004 & 7.3e+00 & 4.8e-09 & 6.8e+01 & 2.3e-08 & F & -- & F & -- & F & -- \\
A2NSDSIL & 60012 & 20004 & F & -- & 4.0e+01 & 1.4e+01 & 1.8e+01 & 1.8e+00 & 1.4e+02 & 1.3e+00 & 6.0e+00 & 5.8e+00 \\
A2NSDSSL & 60012 & 20004 & 6.3e+00 & 6.6e-12 & 3.7e+01 & -7.5e-06 & 5.8e+00 & -6.0e-01 & 1.2e+01 & 6.5e-05 & 1.2e+00 & 4.4e-07 \\
A2NSSSSL & 60012 & 20004 & 3.8e+00 & 0.0e+00 & 2.9e+01 & -1.1e-06 & 2.2e+01 & -1.3e-03 & 1.3e+01 & -1.1e-04 & 9.7e-01 & 5.4e-08 \\
A5ENDNDL & 45006 & 15002 & 6.5e+00 & 1.6e-08 & 4.4e+01 & -7.8e-08 & F & -- & F & -- & F & -- \\
A5ENINDL & 45006 & 15002 & 3.4e+00 & 1.2e-08 & 5.9e+01 & -6.7e-09 & F & -- & F & -- & F & -- \\
A5ENSNDL & 45006 & 15002 & 1.6e+00 & 3.2e-10 & 5.5e+00 & 5.9e-09 & 1.3e+01 & 9.8e-04 & 3.8e+00 & 7.1e-04 & 8.0e-01 & 1.8e-04 \\
A5ESDNDL & 45006 & 15002 & 1.6e+00 & 1.7e-08 & 4.1e+01 & 2.0e-08 & F & -- & F & -- & F & -- \\
A5ESINDL & 45006 & 15002 & 3.1e+00 & 1.0e-08 & 4.3e+01 & -3.9e-08 & F & -- & F & -- & F & -- \\
A5ESSNDL & 45006 & 15002 & 1.2e+00 & 3.1e-08 & 5.1e+00 & -6.1e-09 & 1.4e+01 & 9.7e-04 & 3.3e+00 & 1.2e-03 & 9.7e-01 & 1.8e-04 \\
A5NNDNDL & 60012 & 20004 & 2.8e+02 & 4.2e-06 & 2.1e+03 & 1.7e-06 & F & -- & F & -- & F & -- \\
A5NNDNIL & 60012 & 20004 & 4.9e+00 & 4.1e+01 & 3.4e+01 & 1.0e+01 & F & -- & F & -- & F & -- \\
A5NNDNSL & 60012 & 20004 & 5.1e+01 & 9.9e-08 & 1.6e+02 & -1.1e-07 & 5.9e+00 & -6.2e-01 & 1.3e+01 & -2.2e-04 & 1.8e+00 & 4.3e-05 \\
A5NNSNSL & 60012 & 20004 & 3.1e+00 & 0.0e+00 & 2.9e+01 & -3.5e-07 & 4.2e+01 & -8.8e-01 & 1.2e+01 & 2.0e-04 & 8.6e-01 & 2.9e-05 \\
A5NSDSDL & 60012 & 20004 & 1.2e+01 & 1.5e-08 & 1.1e+02 & 2.1e-08 & F & -- & F & -- & F & -- \\
A5NSDSDM & 6012 & 2004 & 1.1e+00 & 0.0e+00 & 9.4e-01 & -2.4e-07 & 5.3e-01 & -2.6e-05 & 9.6e-01 & 1.9e-05 & 2.7e-01 & 1.1e-11 \\
A5NSDSIL & 60012 & 20004 & 1.8e+01 & 1.8e+01 & 3.3e+01 & 1.4e+01 & 3.0e+01 & 1.6e+00 & 1.1e+02 & 1.0e+00 & 3.2e+00 & 1.8e+00 \\
A5NSDSSL & 60012 & 20004 & 1.1e+01 & 0.0e+00 & 4.4e+01 & -7.0e-07 & 6.4e+00 & -4.4e-03 & 1.5e+01 & -1.1e-03 & 1.0e+00 & 2.0e-05 \\
A5NSSNSM & 6012 & 2004 & 1.1e+00 & 0.0e+00 & 1.0e+00 & -2.4e-07 & 9.3e-01 & -2.6e-05 & 9.4e-01 & 1.9e-05 & 3.5e-01 & 1.1e-11 \\
A5NSSSSL & 60012 & 20004 & 3.6e+00 & 0.0e+00 & 2.8e+01 & -3.6e-06 & 5.6e+01 & -3.9e-03 & 1.3e+01 & -8.4e-04 & 1.3e+00 & 1.5e-06 \\
BIGGSC4 & 4 & 13 & 1.8e-01 & -2.4e+01 & 5.4e-04 & -2.4e+01 & 9.8e-03 & -2.4e+01 & 2.8e+00 & -2.5e+01 & 2.3e-01 & -2.4e+01 \\
BLOCKQP1 & 10010 & 5001 & 1.2e+01 & 5.0e+00 & 1.6e-01 & -5.0e+03 & F & -- & F & -- & 7.0e-01 & -5.0e+03 \\
BLOCKQP2 & 10010 & 5001 & 1.7e+00 & -5.0e+03 & 1.5e-01 & -5.0e+03 & F & -- & 1.7e+00 & -5.0e+03 & 5.3e-01 & -5.0e+03 \\
BLOCKQP3 & 10010 & 5001 & 1.4e+02 & 5.0e+00 & 2.0e+02 & -2.5e+03 & F & -- & F & -- & 1.7e+01 & -2.5e+03 \\
BLOCKQP4 & 10010 & 5001 & 1.3e+00 & -2.5e+03 & 3.5e-01 & -2.5e+03 & F & -- & 1.4e+00 & -2.5e+03 & 5.7e-01 & -2.5e+03 \\
BLOCKQP5 & 10010 & 5001 & 2.5e+02 & 5.0e+00 & 2.6e+02 & -2.5e+03 & F & -- & F & -- & F & -- \\
BLOWEYA & 4002 & 2002 & 2.2e+03 & -4.3e-03 & 1.2e+00 & -1.1e-03 & 5.0e+01 & -2.3e-02 & 2.2e-01 & -2.3e-02 & 3.4e-01 & -2.3e-02 \\
BLOWEYB & 4002 & 2002 & 1.0e+03 & 6.4e-04 & 3.8e-02 & -3.3e-05 & 5.1e+01 & -1.5e-02 & 2.0e-01 & -1.5e-02 & 1.9e-01 & 0.0e+00 \\
BLOWEYC & 22 & 12 & 1.7e+01 & -3.1e+00 & 2.5e-03 & -3.1e+00 & 5.5e-03 & -3.1e+00 & 2.4e-02 & -3.1e+00 & 2.0e-01 & -3.1e+00 \\
CLEUVEN3 & 1200 & 2973 & 2.5e+02 & 1.4e+06 & 3.4e+00 & 4.1e+05 & 2.3e+01 & 2.9e+05 & 1.5e+00 & 2.9e+05 & 2.4e+00 & 2.9e+05 \\
CLEUVEN4 & 1200 & 2973 & 3.2e+00 & 4.1e+06 & 2.6e+03 & 3.8e+05 & 6.0e+01 & 2.9e+05 & 4.7e+00 & 2.9e+05 & 2.7e+01 & 2.9e+05 \\
CLEUVEN5 & 1200 & 2973 & 2.5e+02 & 1.4e+06 & 3.9e+00 & 4.1e+05 & 2.0e+01 & 2.9e+05 & 1.3e+00 & 2.9e+05 & 2.4e+00 & 2.9e+05 \\
CLEUVEN6 & 1200 & 3091 & F & -- & 2.5e+00 & 2.2e+07 & 5.6e+01 & 2.2e+07 & 5.3e-01 & 2.2e+07 & 8.3e+00 & 2.2e+07 \\
FERRISDC & 2200 & 210 & 6.0e-01 & -8.6e-08 & 3.1e+00 & -1.0e-10 & F & -- & F & -- & 2.6e-01 & 0.0e+00 \\
GOULDQP1 & 32 & 17 & 1.4e+00 & -3.5e+03 & 2.1e-03 & -3.5e+03 & 1.8e-02 & -3.5e+03 & 1.2e-01 & -3.5e+03 & 2.2e-01 & -3.5e+03 \\
HATFLDH & 4 & 13 & 1.8e-01 & -2.4e+01 & 5.6e-04 & -2.4e+01 & 8.3e-03 & -2.4e+01 & 9.4e-02 & -2.5e+01 & 2.0e-01 & -2.4e+01 \\
HS44 & 4 & 6 & 5.0e-01 & -1.3e+01 & 6.5e-04 & -1.5e+01 & 1.4e-02 & -1.3e+01 & 5.4e-02 & -1.3e+01 & 1.9e-01 & -1.5e+01 \\
HS44NEW & 4 & 6 & 2.1e-01 & -1.3e+01 & 5.2e-04 & -1.5e+01 & 8.7e-03 & -1.3e+01 & 1.0e-01 & -1.3e+01 & 2.2e-01 & -1.5e+01 \\
LEUVEN2 & 1530 & 2329 & 5.5e+02 & -1.4e+07 & 2.2e+00 & -1.4e+07 & 2.5e+00 & -1.4e+07 & 7.1e-01 & -1.4e+07 & F & -- \\
LEUVEN3 & 1200 & 2973 & F & -- & F & -- & 3.9e+02 & -2.0e+09 & 1.4e+01 & -2.0e+09 & 8.2e+01 & -1.9e+09 \\
LEUVEN4 & 1200 & 2973 & 6.5e+02 & -5.0e+08 & F & -- & F & -- & 1.5e+01 & -2.0e+09 & 4.9e+02 & -2.0e+09 \\
LEUVEN5 & 1200 & 2973 & F & -- & F & -- & 3.8e+02 & -2.0e+09 & 1.4e+01 & -2.0e+09 & 9.5e+01 & -1.9e+09 \\
LEUVEN6 & 1200 & 3091 & F & -- & F & -- & 1.5e+02 & -1.2e+09 & 4.0e+00 & -1.2e+09 & 1.2e+02 & -1.2e+09 \\
LEUVEN7 & 360 & 946 & 2.0e+03 & 6.9e+02 & 5.8e-02 & 6.9e+02 & 3.3e-01 & 6.9e+02 & 3.8e-01 & 6.9e+02 & 3.5e-01 & 6.9e+02 \\
LINCONT & 1257 & 419 & 1.1e-01$^{\mathrm{PI}}$ & -- & 1.0e-01$^{\mathrm{PI}}$ & -- & 9.5e-01$^{\mathrm{PI}}$ & -- & 5.6e-01$^{\mathrm{PI}}$ & -- & 5.9e+00$^{\mathrm{PI}}$ & -- \\
MPC1 & 2550 & 3833 & 8.4e+01 & -2.3e+07 & 2.6e+00 & -2.3e+07 & 6.2e+00 & -2.3e+07 & 2.7e+00 & -2.3e+07 & F & -- \\
MPC2 & 1530 & 2351 & 1.8e+03 & -1.5e+07 & 2.7e+00 & -1.5e+07 & 8.3e-01 & -1.5e+07 & 5.3e-01 & -1.5e+07 & F & -- \\
MPC3 & 1530 & 2351 & 4.7e+02 & -1.5e+07 & 2.3e+00 & -1.5e+07 & 8.4e-01 & -1.5e+07 & 7.9e-01 & -1.5e+07 & F & -- \\
MPC4 & 1530 & 2351 & 1.5e+03 & -1.5e+07 & 3.2e+00 & -1.5e+07 & 7.1e-01 & -1.5e+07 & 3.2e-01 & -1.5e+07 & F & -- \\
MPC5 & 1530 & 2351 & 9.5e+02 & -1.5e+07 & 3.1e+00 & -1.5e+07 & 7.1e-01 & -1.5e+07 & 8.9e-01 & -1.5e+07 & F & -- \\
MPC6 & 1530 & 2351 & 2.0e+02 & -1.5e+07 & 2.9e+00 & -1.5e+07 & 9.9e-01 & -1.5e+07 & 2.6e-01 & -1.5e+07 & F & -- \\
MPC7 & 1530 & 2351 & 1.8e+02 & -1.5e+07 & 2.8e+00 & -1.5e+07 & 1.2e+00 & -1.5e+07 & 7.2e-01 & -1.5e+07 & F & -- \\
MPC8 & 1530 & 2351 & 1.8e+02 & -1.5e+07 & 2.8e+00 & -1.5e+07 & 9.4e-01 & -1.5e+07 & 2.7e-01 & -1.5e+07 & F & -- \\
MPC9 & 1530 & 2351 & 1.6e+02 & -1.5e+07 & 3.1e+00 & -1.5e+07 & 8.8e-01 & -1.5e+07 & 6.2e-01 & -1.5e+07 & F & -- \\
MPC10 & 1530 & 2351 & 3.8e+02 & -1.5e+07 & 3.2e+00 & -1.5e+07 & 1.2e+00 & -1.5e+07 & 3.0e-01 & -1.5e+07 & F & -- \\
MPC11 & 1530 & 2351 & 3.1e+02 & -1.5e+07 & 2.7e+00 & -1.5e+07 & 9.0e-01 & -1.5e+07 & 6.8e-01 & -1.5e+07 & F & -- \\
MPC12 & 1530 & 2351 & 3.0e+02 & -1.5e+07 & 2.8e+00 & -1.5e+07 & 7.7e-01 & -1.5e+07 & 2.6e-01 & -1.5e+07 & F & -- \\
MPC13 & 1530 & 2351 & 2.7e+02 & -1.5e+07 & 2.2e+00 & -1.5e+07 & 8.7e-01 & -1.5e+07 & 7.9e-01 & -1.5e+07 & F & -- \\
MPC14 & 1530 & 2351 & 3.3e+02 & -1.5e+07 & 2.3e+00 & -1.5e+07 & 9.1e-01 & -1.5e+07 & 2.9e-01 & -1.5e+07 & F & -- \\
MPC15 & 1530 & 2351 & 2.7e+02 & -1.5e+07 & 2.4e+00 & -1.5e+07 & 9.6e-01 & -1.5e+07 & 6.0e-01 & -1.5e+07 & F & -- \\
MPC16 & 1530 & 2351 & 2.4e+02 & -1.5e+07 & 2.1e+00 & -1.5e+07 & 9.8e-01 & -1.5e+07 & 3.3e-01 & -1.5e+07 & F & -- \\
NASH & 72 & 24 & 3.5e-01$^{\mathrm{PI}}$ & -- & 8.8e-04$^{\mathrm{PI}}$ & -- & 3.6e-02$^{\mathrm{PI}}$ & -- & 2.5e-01$^{\mathrm{PI}}$ & -- & 0.0e+00$^{\mathrm{PI}}$ & -- \\
NCVXQP1 & 10000 & 5000 & 4.1e+01 & -7.5e+09 & F & -- & F & -- & 1.3e+01 & -7.5e+09 & 2.0e+01 & -7.5e+09 \\
NCVXQP2 & 10000 & 5000 & 1.4e+02 & -5.8e+09 & F & -- & F & -- & 2.6e+01 & -5.8e+09 & 3.0e+03 & -5.8e+09 \\
NCVXQP3 & 10000 & 5000 & 8.2e+01 & -3.1e+09 & F & -- & F & -- & 2.3e+01 & -3.0e+09 & 5.6e+01 & -3.1e+09 \\
NCVXQP4 & 10000 & 2500 & 2.8e+00 & -9.4e+09 & 1.5e+03 & -9.4e+09 & 1.2e+03 & -9.4e+09 & 1.9e+01 & -9.4e+09 & 2.9e+01 & -9.4e+09 \\
NCVXQP5 & 10000 & 2500 & 9.9e+00 & -6.6e+09 & 1.9e+03 & -6.6e+09 & 1.9e+03 & -6.6e+09 & 1.4e+01 & -6.6e+09 & 3.8e+01 & -6.6e+09 \\
NCVXQP6 & 10000 & 2500 & F & -- & 2.7e+03 & -3.4e+09 & 2.8e+03 & -3.4e+09 & 5.3e+01 & -3.4e+09 & 3.0e+01 & -3.5e+09 \\
NCVXQP7 & 10000 & 7500 & 8.0e+01 & -5.2e+09 & F & -- & F & -- & 6.4e+00 & -5.2e+09 & 1.1e+01 & -5.2e+09 \\
NCVXQP8 & 10000 & 7500 & 6.6e+01 & -3.6e+09 & F & -- & F & -- & 2.3e+01 & -3.6e+09 & 2.3e+01 & -3.6e+09 \\
NCVXQP9 & 10000 & 7500 & 1.2e+02 & -2.1e+09 & F & -- & F & -- & 1.4e+01 & -2.1e+09 & 5.7e+01 & -2.1e+09 \\
PORTSNQP & 10 & 2 & 1.5e-01 & 1.9e+00 & 4.2e-04 & 1.9e+00 & 4.8e-03 & 1.9e+00 & 1.0e-01 & 1.9e+00 & 1.5e-01 & 1.9e+00 \\
QPNBAND & 50000 & 25000 & 6.4e-01 & -2.5e+05 & 6.4e+00 & -2.5e+05 & 2.6e+00 & -2.5e+05 & 1.3e-01 & -2.5e+05 & 8.4e-01 & -2.5e+05 \\
QPNBLEND & 83 & 74 & 6.0e-01 & -9.1e-03 & 2.4e-03 & -9.1e-03 & F & -- & F & -- & 1.6e-01 & -9.1e-03 \\
QPNBOEI1 & 384 & 440 & 2.2e+00 & 6.8e+06 & 1.4e-01 & 6.8e+06 & 1.2e+00 & 6.8e+06 & F & -- & 2.3e-01 & 6.8e+06 \\
QPNBOEI2 & 143 & 185 & 2.7e+00 & 1.4e+06 & 1.8e-02 & 1.4e+06 & 4.0e-01 & 1.4e+06 & 2.5e+00 & 1.4e+06 & 3.1e-01 & 1.4e+06 \\
QPNSTAIR & 467 & 356 & 1.4e+00 & 5.1e+06 & 4.4e-02 & 5.1e+06 & 7.1e-01 & 5.1e+06 & 8.5e-01 & 5.1e+06 & F & -- \\
SOSQP1 & 5000 & 2501 & 2.5e-01 & -5.8e-08 & 2.9e-02 & 6.0e-05 & 8.4e-02 & 2.3e-12 & 2.6e-01 & -2.0e-05 & 2.5e-01 & 0.0e+00 \\
SOSQP2 & 5000 & 2501 & 2.3e+01 & -1.2e+03 & F & -- & F & -- & F & -- & 4.7e-01 & -1.2e+03 \\
STATIC3 & 434 & 96 & 4.3e-02$^{\mathrm{DI}}$ & -- & 2.0e-03$^{\mathrm{DI}}$ & -- & F & -- & 3.9e-01$^{\mathrm{DI}}$ & -- & 1.2e-01$^{\mathrm{DI}}$ & -- \\
STNQP1 & 8193 & 4095 & 1.3e+00 & -3.1e+05 & 1.1e+01 & -3.1e+05 & 1.1e+02 & -3.1e+05 & 2.9e-01 & -3.1e+05 & 5.5e-01 & -3.1e+05 \\
STNQP2 & 8193 & 4095 & 1.1e+00 & -5.7e+05 & 3.3e+01 & -5.7e+05 & 5.8e+00 & -5.7e+05 & 7.4e-01 & -5.7e+05 & 6.5e-01 & -5.7e+05 \\
\end{longtable}
\endgroup

\begingroup
\scriptsize
\setlength{\tabcolsep}{1.0pt}
\begin{longtable}{@{}lcc|cc|cc|cc|cc|cc@{}}
\caption{Per-instance million-variable results for $\varepsilon_a=\varepsilon_r=10^{-4}$.  The common start is the frozen $x^0$, and each solver attempt has a $3600$-second limit.}\label{tab:million-detail}\\
\hline
Problem & $n$ & $m$ & \multicolumn{2}{|c}{\PDNQP{}} & \multicolumn{2}{|c}{QPALM} & \multicolumn{2}{|c}{Ipopt} & \multicolumn{2}{|c}{MadNLPGPU} & \multicolumn{2}{|c}{Knitro} \\
\cline{4-13}
& & & Time & Obj & Time & Obj & Time & Obj & Time & Obj & Time & Obj \\
\hline
\endfirsthead
\multicolumn{13}{c}{\tablename\ \thetable\ (continued)} \\
\hline
Problem & $n$ & $m$ & \multicolumn{2}{|c}{\PDNQP{}} & \multicolumn{2}{|c}{QPALM} & \multicolumn{2}{|c}{Ipopt} & \multicolumn{2}{|c}{MadNLPGPU} & \multicolumn{2}{|c}{Knitro} \\
\cline{4-13}
& & & Time & Obj & Time & Obj & Time & Obj & Time & Obj & Time & Obj \\
\hline
\endhead
\hline\multicolumn{13}{r}{Continued on next page}\\
\endfoot
\hline
\endlastfoot
softlabel\_101 & 999999 & 333339 & 3.7e+01 & -1.9e+05 & F & -- & F & -- & F & -- & F & -- \\
softlabel\_102 & 1000000 & 125016 & 1.3e+01 & -8.2e+04 & F & -- & F & -- & F & -- & F & -- \\
softlabel\_103 & 999999 & 333339 & 7.6e+01 & -2.0e+05 & F & -- & F & -- & F & -- & F & -- \\
softlabel\_104 & 1000000 & 200010 & F & -- & F & -- & F & -- & F & -- & F & -- \\
softlabel\_105 & 1000000 & 125016 & 2.8e+01 & -8.9e+04 & F & -- & F & -- & F & -- & F & -- \\
softlabel\_106 & 999999 & 333339 & F & -- & F & -- & F & -- & F & -- & F & -- \\
softlabel\_107 & 1000000 & 100020 & F & -- & F & -- & F & -- & F & -- & F & -- \\
softlabel\_108 & 1000000 & 200010 & 6.9e+01 & -1.2e+05 & F & -- & F & -- & F & -- & F & -- \\
softlabel\_109 & 999999 & 333339 & 4.7e+01 & -2.6e+05 & F & -- & F & -- & F & -- & F & -- \\
softlabel\_110 & 1000000 & 125016 & 3.9e+02 & -7.0e+04 & F & -- & F & -- & F & -- & F & -- \\
flow\_101 & 1000005 & 400002 & 1.2e+02 & 2.6e+05 & F & -- & F & -- & F & -- & F & -- \\
flow\_102 & 999999 & 363636 & 1.4e+02 & 3.0e+05 & F & -- & F & -- & F & -- & F & -- \\
flow\_103 & 999999 & 380952 & 4.5e+01 & 1.8e+05 & F & -- & F & -- & F & -- & F & -- \\
flow\_104 & 999990 & 355552 & 5.8e+01 & 1.5e+05 & F & -- & F & -- & F & -- & F & -- \\
flow\_105 & 999999 & 380952 & 4.1e+02 & 8.7e+04 & F & -- & F & -- & F & -- & F & -- \\
flow\_106 & 999999 & 444444 & 3.4e+01 & 1.7e+05 & F & -- & F & -- & F & -- & F & -- \\
flow\_107 & 999999 & 380952 & 6.9e+01 & 2.6e+05 & F & -- & F & -- & F & -- & F & -- \\
flow\_108 & 1000005 & 400002 & 4.7e+01 & 1.5e+05 & F & -- & F & -- & F & -- & F & -- \\
flow\_109 & 999999 & 444444 & 2.8e+01 & 3.6e+05 & F & -- & F & -- & F & -- & F & -- \\
flow\_110 & 999999 & 444444 & 1.1e+02 & 2.8e+05 & F & -- & F & -- & F & -- & F & -- \\
svm\_101 & 1000000 & 1 & 3.1e+02 & -5.8e+05 & F & -- & F & -- & F & -- & F & -- \\
svm\_102 & 1000000 & 1 & 1.1e+01 & -5.3e+05 & F & -- & F & -- & F & -- & F & -- \\
svm\_103 & 1000000 & 1 & 9.7e+00 & -5.1e+05 & F & -- & F & -- & F & -- & F & -- \\
svm\_104 & 1000000 & 1 & 2.8e+00 & -6.1e+05 & F & -- & F & -- & F & -- & F & -- \\
svm\_105 & 1000000 & 1 & 2.8e+00 & -5.0e+05 & F & -- & F & -- & F & -- & F & -- \\
svm\_106 & 1000000 & 1 & 1.1e+01 & -5.1e+05 & F & -- & F & -- & F & -- & F & -- \\
svm\_107 & 1000000 & 1 & 1.4e+01 & -5.3e+05 & F & -- & F & -- & F & -- & F & -- \\
svm\_108 & 1000000 & 1 & 4.2e+01 & -5.7e+05 & F & -- & F & -- & F & -- & F & -- \\
svm\_109 & 1000000 & 1 & 3.4e+00 & -3.2e+05 & F & -- & F & -- & F & -- & F & -- \\
svm\_110 & 1000000 & 1 & 1.6e+02 & -6.1e+05 & F & -- & F & -- & F & -- & F & -- \\
\end{longtable}
\endgroup
